\documentclass[final,3p]{elsarticle}

\usepackage[colorlinks=true]{hyperref}
\usepackage[hyperref,dvipsnames]{xcolor}
\usepackage[normalem]{ulem}
\usepackage{cancel}
\usepackage[T1]{fontenc} 							
\usepackage[english]{babel} 					
\usepackage{amssymb,amsmath,amsthm,amsbsy}		
\usepackage{array} 						
\usepackage{pdfpages}								
\usepackage{graphicx} 						
\usepackage{caption}[2004/07/16]					
\usepackage{subfigure}							
\usepackage{mathrsfs} 	
\usepackage{empheq}								
\usepackage{nicefrac}								
\usepackage{stmaryrd}								
\usepackage{fixltx2e} 							
\usepackage{fix-cm}								
\usepackage{bm}
\usepackage{booktabs}
\usepackage{here}								
\usepackage{url}
\usepackage{enumitem} 
\usepackage{thmbox}
\usepackage{shadethm}		
\usepackage{siunitx}  		
\usepackage{titlesec}		

\newcommand{\R}{\mathbb{R}}
\newcommand{\N}{\mathbb{N}}
\newcommand\DD[1]{{\mathbb{D}\rb{#1}}}
\newcommand\SP[1]{{\mathbb{S}\rb{#1}}}
\newcommand\Pk[1]{{{\mathbb{P}_{#1}}}}
\newcommand\Qk[1]{{{\mathbb{Q}_{#1}}}}
\newcommand\Pdk[1]{{{\mathbb{P}_{#1}^\mathrm{disc}}}}

\newcommand{\invK}{{\mathrm{inv,K}}}
\newcommand{\inv}{{\mathrm{inv}}}
\newcommand{\maxrm}{{\mathrm{max}}}

\newcommand{\uu}{{\boldsymbol{u}}}
\newcommand{\tuu}{{\tilde{\boldsymbol{u}}}}
\newcommand{\tp}{{\tilde{p}}}
\newcommand{\vv}{{\boldsymbol{v}}}
\newcommand{\vvl}{{\boldsymbol{\hat{v}}_h}}
\newcommand{\ww}{{\boldsymbol{w}}}

\newcommand{\MM}{{\boldsymbol{M}}}

\newcommand{\zero}{{\mathbf{0}}}
\newcommand{\uETA}{{\boldsymbol{\eta}_\uu}}
\newcommand{\uXI}{{\boldsymbol{\xi}_\uu}}
\newcommand{\pETA}{{{\eta}_p}}
\newcommand{\pXI}{{{\xi}_p}}
\newcommand{\uERR}{{\boldsymbol{\mathsf{e}}_\uu}}
\newcommand{\pERR}{{\mathsf{e}_p}}

\newcommand{\ff}{{\boldsymbol{f}}}
\newcommand{\gbld}{{\boldsymbol{g}}}
\newcommand{\x}{{\boldsymbol{x}}}
\newcommand{\n}{{\boldsymbol{n}}}
\newcommand{\TAU}{{\boldsymbol{\tau}}}
\newcommand{\drm}{{\mathrm{d}}}
\newcommand{\Rey}{{\mathrm{Re}}}
\newcommand{\dx}{{\,\drm\x}}
\newcommand{\ds}{{\,\drm\boldsymbol{s}}}

\newcommand{\dtau}{{\,\drm\tau}}
\newcommand{\tend}{{T}}
\newcommand{\half}{{\nicefrac{1}{2}}}
\newcommand{\eps}{{\varepsilon}}
\newcommand{\dvg}{{\mathrm{div}}}

\newcommand{\ih}{{\boldsymbol{i}_h}}
\newcommand{\PI}{{\boldsymbol{\pi}}}
\newcommand{\PIs}{{\boldsymbol{\pi}_s}}

\newcommand{\pis}{{{\pi}_s}}

\newcommand{\fT}{{\mathfrak{T}}}

\newcommand{\ljmp}{\left\llbracket}									
\newcommand{\rjmp}{\right\rrbracket}									
\newcommand\jmp[1]{\ljmp#1\rjmp}										
\newcommand{\Fi}{{\mathcal{F}_h^i}}					

\newcommand\XX{{\boldsymbol{\mathcal{X}}}}

\newcommand\Cc[1]{{\mathcal{C}\rb{#1}}}
\newcommand\CC[1]{{\boldsymbol{\mathcal{C}}\rb{#1}}}
\newcommand\CK[1]{{\boldsymbol{\mathcal{C}}^{#1}}}

\newcommand\Ltwo{{{\mathcal{L}^{2}}}}	
\newcommand\LTWO{{{\boldsymbol{\mathcal{L}}^{2}}}}	
\newcommand\Lp[1]{{{\mathcal{L}^{#1}}}} 
\newcommand\LP[1]{{{\boldsymbol{\mathcal{L}}^{#1}}}} 
\newcommand\Lpz[1]{{{\mathcal{L}_0^{#1}}}} 

\newcommand\Wmp[2]{{{\mathcal{W}^{#1,#2}}}}

\newcommand\WMP[2]{{{\boldsymbol{\mathcal{W}}^{#1,#2}}}}

\newcommand\Hm[1]{{{\mathcal{H}^{#1}}}}
\newcommand\Hmz[1]{{{\mathcal{H}_0^{#1}}}}
\newcommand\HM[1]{{{\boldsymbol{\mathcal{H}}^{#1}}}}
\newcommand\HMZ[1]{{{\boldsymbol{\mathcal{H}}_0^{#1}}}}

\newcommand\HDIV[1]{{{\boldsymbol{\mathcal{H}}\rb{\dvg;#1}}}}
			
\newcommand{\VV}{{\boldsymbol{\mathcal{V}}}}							
\newcommand{\Q}{{\mathcal{Q}}}
\newcommand{\T}{{\mathcal{T}_h}} 

\newcommand\rb[1]{\hspace{-0.15ex}\left(#1\right)}
\newcommand\ger[1]{\mathrm{#1}}
\newcommand\sqb[1]{\left[ #1 \right]}

\newcommand\set[1]{\left\{ #1 \right\}}
\newcommand\bra[1]{\langle #1 \rangle}
\newcommand\abs[1]{\left\lvert#1\right\rvert}
\newcommand\norm[1]{\left\lVert#1\right\rVert}
\newcommand\enorm[1]{\norm{#1}_e}

\newcommand{\tripnorm}[1]{{\left\vert\kern-\nulldelimiterspace\left\vert\kern-\nulldelimiterspace\left\vert #1
	\right\vert\kern-\nulldelimiterspace\right\vert\kern-\nulldelimiterspace\right\vert}}
	
\newcommand{\otoprule}{\midrule[\heavyrulewidth]}

\newcommand\restr[2]{{												
	\left.\kern-\nulldelimiterspace									
	#1
	\vphantom{\big|}
	\right|_{#2}
	}}
	
\newcommand{\goodgap}{%
	\hspace{0.01\subfigtopskip}
	\hspace{0.01\subfigbottomskip}
	}							
	
\newtheorem[style=S,underline=true,bodystyle=\normalfont\noindent]{thmDef}{\textsc{Definition}}[section]
\newtheorem[style=S,cut=false]{thmCor}[thmDef]{\textsc{Corollary}}
\newtheorem[style=S,cut=false,headstyle=\normalsize\bfseries\boldmath####1~####2]{thmLem}[thmDef]{\textsc{Lemma}}

\newshadetheorem{thm}[thmDef]{\textsc{Theorem}}
\newenvironment{thmThe}[1][]{%
 	\definecolor{shadethmcolor}{gray}{0.9}%
  	\definecolor{shaderulecolor}{rgb}{0.0,0.0,0.0}%
  	\setlength{\shadeboxrule}{0.0pt}%
  	\begin{thm}[#1]%
}{\end{thm}}			

\newenvironment{thmProof}
                [0]
                { \begin{example}[\textsc{Proof}] \normalsize}
                { $\hfill\blacksquare$ \end{example} }

\newenvironment{thmRem}
                [0]
                { \refstepcounter{thmDef} \begin{example}[\textsc{Remark} \thesection.\arabic{thmDef}]  \normalsize}
                { $\hfill\blacktriangle$ \end{example} }		
                
\titleformat{\section}{\fontsize{10.5}{17}\bfseries}{\thesection}{1em}{}
\titleformat{\subsection}{\bfseries\itshape}{\thesubsection}{1em}{}

\journal{the `Journal of Numerical Mathematics' (accepted: January 29, 2017)}

\definecolor{mediumblue}{RGB}{0,0,205}
\definecolor{forestgreen}{RGB}{34,139,34}
\definecolor{darkred}{RGB}{200,0,0}

\begin{document}

\hypersetup{
  linkcolor=darkred,
  urlcolor=forestgreen,
  citecolor=mediumblue
}


\begin{frontmatter}


\title{Pressure-robust analysis of divergence-free and conforming FEM \\ for evolutionary incompressible Navier--Stokes flows}

\author[]{Philipp W.\ Schroeder\corref{cor1} \fnref{fn1}}
	\cortext[cor1]{Corresponding author}
	\fntext[fn1]{ORCID: \url{https://orcid.org/0000-0001-7644-4693}}		
    \ead{p.schroeder@math.uni-goettingen.de}

\author[]{Gert Lube} 
	\ead{lube@math.uni-goettingen.de}

\address{Institute for Numerical and Applied Mathematics, Georg-August-University G\"ottingen, D-37083 G\"ottingen, Germany}

\begin{abstract}

This article focusses on the analysis of a conforming finite element method for the time-dependent incompressible Navier--Stokes equations. For divergence-free approximations, in a semi-discrete formulation, we prove error estimates for the velocity that hold independently of both pressure and Reynolds number. Here, a key aspect is the use of the discrete Stokes projection for the error splitting. Optionally, edge-stabilisation can be included in the case of dominant convection. Emphasising the importance of conservation properties, the theoretical results are complemented with numerical simulations of vortex dynamics and laminar boundary layer flows.

\vspace{0.25cm}
 
\end{abstract}

\begin{keyword}

Incompressible viscous flow \sep
divergence-free FEM \sep
pressure/semi-robust error estimates \sep
vortex dynamics

\vspace{0.25cm}
\noindent \emph{Short title:} Pressure-robust analysis of div-free conforming FEM for incompressible flows

\end{keyword}

Publisher's version: DOI \url{https://doi.org/10.1515/jnma-2016-1101}

\end{frontmatter}

\section{Introduction}	\label{sec:Intro}

In this paper, we consider the time-dependent incompressible Navier--Stokes equations (NSEs) \cite{Tritton88,SchlichtingGersten00,Durst08}
\begin{subequations}\label{eq:TINS}
	\begin{empheq}[left=\empheqlbrace]{alignat=2} 
		\partial_t\uu - \nu\Delta \uu + \rb{\uu\cdot\nabla}\uu +\nabla p &= \ff \qquad\quad 												&&\text{in }\rb{0,\tend}\times\Omega, 			\\
		\nabla\cdot\uu &= 0 				&&\text{in }\rb{0,\tend}\times\Omega, 			\\
		\uu &= \zero 				&&\text{on }\rb{0,\tend}\times\partial\Omega, \\
		\uu\rb{0,\x} &=\uu_0\rb{\x} 	&&\text{for } \x\in\Omega,		
	\end{empheq} 
\end{subequations}
with no-slip condition in a bounded, polyhedral and convex Lipschitz domain $\Omega\subset\R^d$ for the space dimension $d\in\set{2,3}$. Here, $\uu \colon\rb{0,\tend}\times\Omega\to\R^d$ denotes the velocity field, $p\colon\rb{0,\tend}\times\Omega\to\R$ is the (zero-mean) kinematic pressure, $\ff\colon\rb{0,\tend}\times\Omega\to\R^d$ represents external body forces and $\uu_0\colon \Omega\to\R^d$ stands for a suitable initial condition for the velocity. The underlying fluid is assumed to be a Newtonian fluid with constant dimensionless kinematic viscosity $\nu\geqslant 0$, where $\nu=0$ corresponds to the incompressible Euler equations.\\

In this work, we consider the spatial approximation of \eqref{eq:TINS} by `weakly divergence-free', $\HM{1}$-conforming and inf-sup stable finite element methods (FEMs). If $\VV_h$ and $\Q_h$ denote the finite element spaces for the discrete velocity $\uu_h$ and pressure $p_h$, respectively, weakly divergence-free FEM are characterised by the inclusion property $\nabla\cdot\VV_h\subset\Q_h$, which leads to $\nabla\cdot\uu_h=0$ pointwise. A major advantage of such methods, for example shown in \cite{JohnEtAl16} for the Stokes problem, is that the pressure approximation does not influence the velocity approximation; a property called `pressure-robustness'. Concerning the regularity assumptions, we require $\uu\in\Lp{1}{\rb{0,\tend;\WMP{1}{\infty}}}$ which is a rather natural assumption in the context of FEM; cf. \cite{MajdaBertozzi02,BurmanFernandez07,GiraultEtAl15}.\\

The first aim of this paper is to perform a detailed semi-discrete numerical analysis, including stability analysis and error estimates, for the continuous-in-time solution of weakly divergence-free, conforming FEM. In doing so, we attach great importance to ensuring that the constants in our estimates do not explicitly depend on $\nu^{-1}$. That is, provided the exact solution is sufficiently smooth, the analysis allows for estimates which hold uniformly in the Reynolds number and therefore also for the incompressible Euler equations. Dependence on $\nu^{-1}$ can only be seen implicitly through certain Sobolev norms of the exact solution and, inspired by \cite{RoosEtAl08}, we call such estimates `semi-robust'. \\

Let us give a short overview of some previous research in this direction and highlight similarities and differences. Note that the following selection is restricted to $\HM{1}$-conforming methods.

\begin{itemize}
\item  \textbf{Steady Oseen problem:} An equal-order, continuous interior penalty (CIP) method with edge-based stabilisation of the jumps of both the gradient in normal direction and the divergence is presented in \cite{BurmanEtAl06}, where the first stabilisation accounts for dominant convection and the second one gives additional control of the incompressibility constraint. In \cite{MatthiesTobiska15}, inf-sup stable FEMs with local projection stabilisation (LPS) for controlling convection as well as the divergence of the discrete velocity are considered. Furthermore, the CIP-based method penalises jumps of the pressure gradient to become stable whereas the LPS-based method, for discontinuous pressures, optionally includes a pressure jump penalisation term to recover $\mathcal{O}\rb{h^\half}$ in the convergence analysis. In \cite{BurmanLinke08} a Scott--Vogelius (SV) FEM with either LPS or edge-based stabilisation is analysed. SV-FEM fulfil $\nabla\cdot\VV_h\subset\Q_h$ and, therefore, divergence stabilisation is redundant. On special macro triangulations, SV-FEM are inf-sup stable and no pressure stabilisation is necessary. Actually, any kind of pressure stabilisation would modify the discrete mass balance and in doing so destroy the exact mass conservation on the discrete level.
\item \textbf{Evolutionary Oseen problem:} Major contributions to the analysis of equal-order methods based on orthogonal subscales can be found in \cite{Codina02} and the references therein. For inf-sup stable FE pairs, \cite{FrutosEtAl16} considers grad-div stabilisation exclusively and in addition to semi-discrete analysis also provides error estimates for commonly used time-stepping schemes (implicit Euler, BDF(2) and Crank--Nicolson). An LPS method with stabilisation for both incompressibility constraint and dominant convection is analysed in \cite{DallmannEtAl16}, wherein some results from \cite{MatthiesTobiska15} are extended and improved.
\item \textbf{Evolutionary Navier--Stokes problem:} Previous research where the full problem \eqref{eq:TINS} is considered from the point of view of semi-robustness is presented in \cite{BurmanFernandez07,ArndtEtAl15}. Being consequent extensions of \cite{BurmanEtAl06,DallmannEtAl16}, respectively, the former work considers edge-stabilisation for equal-order methods, whereas the latter relies on LPS for stabilising divergence, dominant convection and, for discontinuous pressure approximations, pressure jumps. Assuming that both exact velocity and pressure are spatially smooth according to $\uu,p\in\HM{k+1}{\rb{\Omega}}$, the $\LTWO$-velocity error of the edge-stabilised equal-order method \cite{BurmanFernandez07} is shown to converge quasi-optimally with $\mathcal{O}\rb{h^{k+\half}}$. Note that, due to using the $\LTWO$-projection as an approximation operator, no assumption for the regularity of $\partial_t\uu$ is needed. The analysis for the inf-sup stable LPS method \cite{ArndtEtAl15}, however, is based on an approximation operator that preserves the discrete divergence and, in the standard case, has a suboptimal convergence rate $\mathcal{O}\rb{h^{k}}$, provided that $\uu\in\HM{k+1}{\rb{\Omega}}$ and $\partial_t\uu,p\in\HM{k}{\rb{\Omega}}$. For high Reynolds numbers, similarly to the steady Oseen problem, it is possible to recover $\mathcal{O}\rb{h^\half}$ in the convergence analysis by means of an enriched discrete pressure space, whenever additionally $\partial_t\uu\in\HM{k+1}{\rb{\Omega}}$. Moreover, for the two-dimensional Navier--Stokes problem, \cite{Burman15} shows that it is possible to use the stream function/vorticity formulation to obtain semi-robust FEM error estimates. Interestingly, the analysis is based on a scale separation into large eddies and small scales known from LES modelling.

\end{itemize}

Our second aim is to augment the theoretical results by a thorough analysis of some selected numerical examples. In the literature, several CFD benchmarks can be found which demonstrate that the lack of pressure-robustness of more popular numerical methods can cause severe practical problems. Exemplarily, we want to refer to \cite{JohnEtAl16} where the no-flow problem, a stationary vortex, a flow with Coriolis force and natural convection has been shown to benefit from pressure-robust methods. Convincing results in this direction are also reported in \cite{LinkeMerdon16,LedererEtAl16} where pressure-robust reconstruction techniques are applied advantageously to potential flows. In addition, emphasising the importance of the treatment of the inertia term, in \cite{CharnyiEtAl17} a new formulation with improved conservation properties for the nonlinear term has been developed. Numerical tests with the Gresho-vortex problem and several flows over immersed bodies indicate that this aspect of discretisation schemes is also important for the successful application of a numerical method. Our contribution lies in comparing weakly divergence-free methods with several versions of the standard Taylor--Hood element. We find that while in the case of flows with vortical structure the divergence-free methods are clearly superior, for the laminar Blasius boundary layer they are at least as good as the Taylor--Hood type methods.  \\

\textbf{Organisation of the article:} In Section \ref{sec:SpecialClass}, we introduce weakly divergence-free, $\HM{1}$-conforming and inf-sup stable FEM for the time-dependent Navier--Stokes problem. Then, in Section \ref{sec:StabWellPosed}, we briefly discuss the stability and well-posedness of the method. The main part of this work, Section \ref{sec:ErrorEstimatesVel}, is concerned with the derivation of a-priori, pressure- and semi-robust error estimates for the discrete velocity. Some numerical results for viscous vortex flows and laminar boundary layer flows are presented in Section \ref{sec:NumEx}. Supplementary, an error analysis for the corresponding pressure approximation can be found in Section \ref{sec:ErrorEstimatesPre}. 

\section{A special class of FEM for incompressible flows} \label{sec:SpecialClass}	
In this section, we recall the weak formulation with its functional analytic background and give regularity assumptions for the incompressible flow problem \eqref{eq:TINS}. Afterwards, we introduce the key aspects of the space semi-discretisation of this problem using $\HM{1}$-conforming, weakly divergence-free and inf-sup stable FEM.

\subsection{Incompressible Navier--Stokes equations}	
In what follows, for $K\subseteq\Omega$, we frequently use the standard Sobolev space $\Wmp{m}{p}{\rb{K}}$ for scalar-valued functions with associated norm $\norm{\cdot}_{\Wmp{m}{p}{\rb{K}}}$ and seminorm $\abs{\cdot}_{\Wmp{m}{p}{\rb{K}}}$ for $m\geqslant 0$ and $p\geqslant 1$. Spaces and norms for vector- and matrix-valued functions are indicated with bold letters. In particular, we obtain the Lebesgue space $\Wmp{0}{p}{\rb{K}}=\Lp{p}{\rb{K}}$ and the Hilbert space $\Wmp{m}{2}{\rb{K}}=\Hm{m}{\rb{K}}$. Additionally, the closed subspaces $\Hmz{1}{\rb{K}}$ consisting of $\Hm{1}{\rb{K}}$-functions with vanishing trace on $\partial K$ and the set  $\Lpz{2}{\rb{K}}$ of $\Lp{2}{\rb{K}}$-functions with zero mean in $K$ play an important role. The $\Lp{2}{\rb{K}}$-inner product is denoted by $\rb{\cdot,\cdot}_K$ and, for brevity when no confusion can arise, if $K=\Omega$ we usually omit the domain completely. Furthermore, given a Banach space $\XX$, the Bochner space $\Lp{p}{\rb{0,\tend;\XX}}$ for $p\in\sqb{1,\infty}$ is used. With the obvious modification for $p=\infty$, and for $\tripnorm{\cdot}_\XX$ denoting either the norm $\norm{\cdot}_\XX$ or the seminorm $\abs{\cdot}_\XX$, we define
\begin{align}
\tripnorm{\vv}_{\Lp{p}{\rb{0,\tend;\XX}}}^p
	= \int_0^\tend \tripnorm{\vv\rb{\tau}}_\XX^p \dtau.
\end{align}

Now, we can introduce suitable functions spaces for velocity and pressure, respectively, by
\begin{equation}
	\VV=\HMZ{1}{\rb{\Omega}} \qquad \text{and} \qquad
	\Q=\Lpz{2}{\rb{\Omega}}
\end{equation}
and obtain the following continuous variational formulation of problem \eqref{eq:TINS}: 
\begin{subequations}\label{eq:VarTINS}
	\begin{empheq}[left=\empheqlbrace]{align} 
	\text{find }\rb{\uu,p}\colon\rb{0,\tend}\to\VV\times\Q
		\text{ with }\uu\rb{0} =\uu_0	&\text{ s.t. }\forall\rb{\vv,q}\in\VV\times\Q	\\
		\bra{\partial_t\uu,\vv} + a\rb{\uu,\vv} + t\rb{\uu;\uu,\vv} + b\rb{\vv,p} &= \bra{\ff,\vv}  	\\
		  -b\rb{\uu,q} &= 0			
	\end{empheq} 
\end{subequations}
Here, $\bra{\cdot,\cdot}$ denotes the duality pairing between $\VV$ and its dual space $\VV^*=\HM{-1}{\rb{\Omega}}$. The corresponding multilinear forms are given by
\begin{align}
	a\rb{\uu,\vv} = \int_\Omega \nu\nabla\uu:\nabla\vv\dx,\quad
	t\rb{\ww;\uu,\vv} = \int_\Omega \rb{\ww\cdot\nabla}\uu\cdot\vv\dx,\quad
	b\rb{\vv,q}=-\int_\Omega q\rb{\nabla\cdot\vv}\dx.
\end{align}
In the following, we recall some well-known properties of the trilinear form.
\begin{thmLem}[The trilinear form] \label{lem:Trilinear}
Let $\uu,\vv,\ww\in\VV$ with $\nabla\cdot\uu=0$ and $1\leqslant p,q,r\leqslant\infty$ with $\nicefrac{1}{p}+\nicefrac{1}{q}+\nicefrac{1}{r}=1$. Then, the trilinear form is continuous on $\VV\times\VV\times\VV$ and
\begin{subequations}
\begin{align}
t\rb{\uu;\vv,\ww}&=-t\rb{\uu;\ww,\vv},\quad
t\rb{\uu;\vv,\vv}=0, \label{eq:TrilinearSkewSymm}\\
\abs{t\rb{\uu;\vv,\ww}}&\leqslant \norm{\uu}_\LP{p}\norm{\nabla\vv}_\LP{q}\norm{\ww}_\LP{r}.
\end{align}		
\end{subequations}
\end{thmLem}
\begin{thmProof}
Cf., for example, \cite[Lemma V.1.1]{BoyerFabrie13} and \cite[Section 6.2, Lemma 13]{Layton08}.	
\end{thmProof}
The divergence constraint in \eqref{eq:VarTINS} prompts us to define the subspace
\begin{equation}
	\VV^\dvg = \set{\vv\in\VV\colon -b\rb{\vv,q}=\rb{q,\nabla\cdot\vv}=0~\forall q\in\Q} \subset \VV
\end{equation}
of weakly divergence-free functions. Lastly, concerning the regularity of both data and solution, we assume
\begin{align}
\ff\in\Lp{2}\rb{0,\tend;\VV^*}, \quad
\uu\in\Lp{1}\rb{0,\tend;\WMP{1}{\infty}}, \quad
\partial_t\uu\in\Lp{2}{\rb{0,\tend;\LP{2}}}	 \quad \text{and} \quad
\uu_0\in\LP{2},	
\end{align}
which ensures uniqueness of the continuous weak solution as shown in the following lemma. As a motivation for this, let us mention that, for the stationary case, \cite{GiraultEtAl15} discusses bounds of $\abs{\uu}_\WMP{1}{\infty}$ against the data. An analogous extension to the time-dependent case is possible, albeit implying considerable work.

\begin{thmLem}[Uniqueness]
If a solution $\uu\in\Lp{1}{\rb{0,\tend;\WMP{1}{\infty}}}\cap\Hm{1}{\rb{0,\tend;\LP{2}}}$ to the NSEs \eqref{eq:VarTINS} exists, it is unique.
\end{thmLem}

\begin{thmProof}
This proof is based on \cite[Section 3.1]{MajdaBertozzi02}, where only the Cauchy problem is considered. Let $\rb{\uu_1,p_1}$ and $\rb{\uu_2,p_2}$ be two solutions of \eqref{eq:VarTINS} with smooth velocities according to $\Lp{1}{\rb{0,\tend;\WMP{1}{\infty}}}\cap\Hm{1}{\rb{0,\tend;\LP{2}}}$ and initial values $\uu_{0_1}$ and $\uu_{0_2}$. Denote the difference of the solutions by $\rb{\tuu,\tp}=\rb{\uu_1-\uu_2,p_1-p_2}$.\\

Considering the difference of the PDEs \eqref{eq:TINS} for $\rb{\uu_1,p_1}$ and $\rb{\uu_2,p_2}$, respectively, and adding a zero yields
\begin{align}
	\partial_t\tuu + \rb{\uu_1\cdot\nabla}\tuu+ \rb{\tuu\cdot\nabla}\uu_2
		=-\nabla \tp + \nu\Delta \tuu. 
\end{align}

Multiplication by $\tuu$ and integration over $\Omega$ yields
\begin{align}
	\rb{\partial_t\tuu,\tuu} + t\rb{\uu_1;\tuu,\tuu}+ t\rb{\tuu;\uu_2,\tuu}
		=-\rb{\nabla \tp,\tuu} + \nu\rb{\Delta \tuu,\tuu}, 
\end{align}

where $t\rb{\uu_1;\tuu,\tuu}=0$ due to Lemma \ref{lem:Trilinear}. The two terms on the right-hand side can be treated using integration by parts, where the boundary terms cancel out for no-slip or periodic boundary conditions on $\partial\Omega$. With the estimate in Lemma \ref{lem:Trilinear} we thus obtain
\begin{equation}
	\frac{1}{2}\frac{\drm}{\drm t}\norm{\tuu}_\LP{2}^2 
	\leqslant \frac{1}{2}\frac{\drm}{\drm t}\norm{\tuu}_\LP{2}^2 
		+ \nu\norm{\nabla\tuu}_\LP{2}^2 
		= -t\rb{\tuu;\uu_2,\tuu}
		\leqslant \norm{\nabla\uu_2}_\LP{\infty} \norm{\tuu}_\LP{2}^2.
\end{equation}

Applying the differential form of Gronwall's lemma gives, for all $0\leqslant t\leqslant T$,
\begin{equation}
	\norm{\tuu\rb{t}}_\LP{2}^2
		\leqslant \norm{\tuu\rb{0}}_\LP{2}^2
		\exp\rb{2\int_0^t \norm{\nabla\uu_2\rb{\tau}}_\LP{\infty} \dtau},
\end{equation}

which, after taking the square root, becomes
\begin{equation}
	\norm{\tuu}_{\Lp{\infty}\rb{0,\tend;\LP{2}}}
		\leqslant \norm{\tuu\rb{0}}_\LP{2}
		e^{\abs{\uu_2}_{\Lp{1}\rb{0,\tend;\WMP{1}{\infty}}}}.
\end{equation}

The claim follows immediately whenever $\uu_{0_1}$ and $\uu_{0_2}$ coincide.

\end{thmProof}

\subsection{Weakly divergence-free, conforming and inf-sup stable FEM}	
In this work, our aim is to consider conforming finite element approximations of the variational formulation \eqref{eq:VarTINS} by introducing finite-dimensional spaces $\VV_h\subset\VV$ and $\Q_h\subset\Q$ of piecewise degree $k\geqslant 1$ and $\ell\geqslant 0$ polynomials for the discrete velocity and pressure, respectively. The subscript $h$ refers to a quasi-uniform and exact decomposition $\T$ of the domain $\Omega$ without hanging nodes. Also, $h=\max_{K\in\T}h_K$ where $h_K$ denotes the diameter of the particular element $K\in\T$.  It is well-known that for shape-regular decompositions $\T$, the discrete space $\VV_h$ satisfies the local inverse inequality \cite[Lemma 1.138]{ErnGuermond04} 
\begin{equation} \label{eq:LocInvEq}
\forall\vv_h\in\VV_h\colon\quad
\norm{\vv_h}_{\WMP{\ell}{p}{\rb{K}}}
	\leqslant C_\invK h_K^{m-\ell+d\rb{\frac{1}{p}-\frac{1}{q}}}\norm{\vv_h}_{\WMP{m}{q}{\rb{K}}}, \quad
\forall K\in\T,
\end{equation}
where $0\leqslant m\leqslant \ell$ and $1\leqslant p,q\leqslant\infty$. Furthermore, we define $C_\inv=\max_{K\in\T} C_\invK$. Here, the decomposition $\T$ can either consist of simplices or tensor-product elements and we assume that the finite element spaces possess the following optimal approximation properties. There is a velocity approximation operator $\ih\colon \VV\to\VV_h$ such that for all $\vv\in\VV\cap\HM{r}{\rb{\Omega}}$ with $r\geqslant 2$ and $r_\uu=\min\set{r,k+1}$
\begin{align} \label{eq:uInt}
\norm{\vv-\ih\vv}_{\LP{2}{\rb{K}}}+h_K\abs{\vv-\ih\vv}_{\HM{1}{\rb{K}}}
	\leqslant Ch_K^{r_\uu}\abs{\vv}_{\HM{r_\uu}{\rb{K}}},\quad \forall K\in\T.	
\end{align}
Furthermore, concerning the pressure, it is well-known that for all $q\in\Q\cap\HM{s}{\rb{\Omega}}$ with $s\geqslant 1$ and $r_p=\min\set{s,\ell+1}$ the orthogonal $\Ltwo$-projection $\pi_0\colon \Q\to\Q_h$ onto the discrete pressure space fulfils 
\begin{align} \label{eq:pProj}
\norm{q-\pi_0 q}_{\Lp{2}{\rb{K}}}
	\leqslant Ch_K^{r_p}\abs{q}_{\Hm{r_p}{\rb{K}}},\quad \forall K\in\T.
\end{align}
Furthermore, and very importantly, we assume that the finite element spaces fulfil the inclusion property
\begin{equation} \label{eq:divVhinQh}
	\nabla\cdot\VV_h\subseteq\Q_h.
\end{equation}	

Then, the space-semidiscrete variational formulation of \eqref{eq:VarTINS} reads as follows: 
\begin{subequations}\label{eq:DiscVarTINS}
	\begin{empheq}[left=\empheqlbrace]{align} 
	\text{find }\rb{\uu_h,p_h}\colon\rb{0,\tend}\to\VV_h\times\Q_h
	\text{ with }\uu_h\rb{0} =\uu_{0h}&\text{ s.t. }\forall\rb{\vv_h,q_h}\in\VV_h\times\Q_h\\
	\rb{\partial_t\uu_h,\vv_h} + a\rb{\uu_h,\vv_h}+t\rb{\uu_h;\uu_h,\vv_h}+b\rb{\vv_h,p_h}&=\bra{\ff,\vv_h}\\
	-b\rb{\uu_h,q_h} &= 0 				
	\end{empheq} 
\end{subequations}
Here, $\uu_{0h}$ denotes an approximation of $\uu_0$ belonging to $\VV_h$, for example the Stokes projection \cite{GiraultEtAl05}. \\

Introducing the space of discretely divergence-free functions
\begin{equation}
	\VV_h^\dvg = \set{\vv_h\in\VV_h\colon -b\rb{\vv_h,q_h}=\rb{\nabla\cdot\vv_h,q_h}=0~\forall q_h\in\Q_h}	,
\end{equation}
we note that the solution $\uu_h$ of \eqref{eq:DiscVarTINS} is by construction discretely divergence-free. Due to \eqref{eq:divVhinQh} we infer that $\nabla\cdot\uu_h\in\Q_h$. Thus, taking $\nabla\cdot\uu_h\in\Q_h$ as a discrete pressure test function, we deduce that $\uu_h$ is divergence-free in the $\Ltwo$-sense. Moreover, since the velocity of an $\HM{1}$-conforming approximation is globally continuous, that is $\VV_h\subset\CC{\overline{\Omega}}$, we know that $\uu_h$ is divergence-free even pointwise and hence
\begin{equation}
	\VV_h^\dvg = \set{\vv_h\in\VV_h\colon \nabla\cdot\vv_h\rb{\x}=0~\forall\x\in\Omega}	\subset\VV^\dvg.
\end{equation}

\begin{thmRem} \label{rem:ConsTrilinearForm}
Note that the trilinear form $t\rb{\cdot;\cdot,\cdot}$ is used for both the continuous problem \eqref{eq:VarTINS} and the semi-discrete problem \eqref{eq:DiscVarTINS}. Contrary to non-divergence-free FEM, it is not necessary to modify the trilinear form because skew-symmetry and  important conservation properties \cite{CharnyiEtAl17}, explained in more detail in Section \ref{sec:ConservationProperties}, hold automatically whenever the first argument is divergence-free; cf.\ Lemma \ref{lem:Trilinear}. However, the assembly process is computationally more efficient when using the non-modified trilinear form \cite[Section 4.1]{Gunzburger89}, and the analysis simplifies considerably at various points. 
\end{thmRem} 

The last ingredient of the considered class of FEM is inf-sup stability. This means that we assume that the spaces $\VV_h$ and $\Q_h$ satisfy a discrete inf-sup condition
\begin{equation} \label{eq:InfSup}
	\inf_{q_h\in\Q_h\backslash\set{0}} \sup_{\vv_h\in\VV_h\backslash\set{\zero}}
	\frac{\rb{q_h,\nabla\cdot\vv_h}}{\norm{\nabla\vv_h}_\LP{2}\norm{q_h}_\Lp{2}}
	\geqslant \beta_h	,
\end{equation}
for a constant $\beta_h>0$ independent of $h$. Moreover, due to the \eqref{eq:InfSup} and the closed range theorem, $\VV_h^\dvg\neq\set{\zero}$ and we thus avoid locking phenomena \cite{GiraultRaviart86}.

\begin{thmRem}
Note that the gradient of a velocity field is actually a second order tensor and in writing \eqref{eq:InfSup} we implicitly used the equivalence of the full norm $\norm{\vv}_\Hm{1}$ and the seminorm $\abs{\vv}_\HM{1}=\norm{\nabla \vv}_\LP{2}$ on $\HMZ{1}{\rb{\Omega}}$, due to Poincare--Friedrich's inequality \cite[Proposition III.2.38]{BoyerFabrie13}.	
\end{thmRem}

\begin{thmRem}
One example of finite element spaces which fulfil the above conditions on simplicial decompositions is the Scott--Vogelius (SV) pair \cite{ScottVogelius85} of order $k\in\N$. With 
\begin{subequations}
\begin{align}
\Pk{k}&=\set{v_h\in\Cc{\overline{\Omega}}\colon \restr{v_h}{K}\in\Pk{k}\rb{\overline{K}},~\forall K\in\T}, \\
\Pdk{k-1}&=\set{q_h\in\Lp{2}{\rb{\Omega}}\colon \restr{q_h}{K}\in\Pk{k-1}\rb{K},~\forall K\in\T},
\end{align}	
\end{subequations}
for the velocity components and pressure, respectively, we obtain the discrete function spaces by intersection:
\begin{equation}
\VV_h=\sqb{\Pk{k}}^d\cap\VV=\sqb{\Pk{k}}^d\cap \HMZ{1}{\rb{\Omega}}, \qquad
\Q_h=\Pdk{k-1} \cap \Q	=\Pdk{k-1} \cap \Lpz{2}{\rb{\Omega}}	
\end{equation}
Inf-sup stability of SV elements is guaranteed on meshes without singular vertices and polynomial degree $k\geqslant 4$ ($d=2$), for $k\geqslant 6$ ($d=3$) and on barycentre-refined meshes for $k\geqslant d$; cf. \cite{Qin94,Zhang05,Zhang11}. Note that a once barycentre-refined mesh fulfils the local inverse inequality \eqref{eq:LocInvEq} even though the constant might be very large. We will use this element for numerical simulations in Section \ref{sec:NumEx}.
\end{thmRem}

\begin{thmRem}
Of course, there are other spaces which fulfil the above conditions. In \cite{GuzmanNeilan14a,GuzmanNeilan14b}, conforming and divergence-free Stokes elements for $d\in\set{2,3}$ on simplicial meshes are introduced, whose construction is based on enriching $\HDIV{\Omega}$-conforming FE spaces with divergence-free rational functions that enforce strong tangential continuity. On tensor-product meshes for $d=2$, \cite{Zhang09,HuangZhang11} shows that the conforming element $\VV_h/\Q_h=\rb{\Qk{k,k-1}\times\Qk{k-1,k}}/\nabla\cdot\VV_h$ is inf-sup stable. Pointwise divergence-free velocity approximations are ensured by polynomials of different degrees with respect to the coordinate directions, where $\Qk{k,k-1}$ stands for the space of continuous piecewise polynomials of degree $k$ in $x_1$-direction and degree $k-1$ in $x_2$-direction. In the context of IsoGeometric analysis, several conforming and divergence-free finite element spaces have been constructed using splines on tensor-product meshes \cite{BuffaEtAl11,EvansHughes13}.	
\end{thmRem}

Lastly, an important result is stated which will be used frequently for the error estimates.
\begin{thmCor}[Galerkin orthogonality] \label{cor:GalOrtho}
	Let $\rb{\uu,p}\in\VV\times\Q$ solve \eqref{eq:VarTINS} and $\rb{\uu_h,p_h}\in\VV_h\times Q_h$ solve \eqref{eq:DiscVarTINS}. Then, with $\pXI=p-p_h$,
	\begin{align} \label{eq:GalOrtho}
	0=\rb{\partial_t\sqb{\uu-\uu_h},\vv_h} &+ a\rb{\uu-\uu_h,\vv_h} +t\rb{\uu;\uu,\vv_h}
	-t\rb{\uu_h;\uu_h,\vv_h} +b\rb{\vv_h,\pXI} -b\rb{\uu-\uu_h,q_h}
	\end{align}		
	holds almost everywhere in $\rb{0,\tend}$ and for all $\rb{\vv_h,q_h}\in\VV_h\times Q_h$.
\end{thmCor}
\begin{thmProof}
Subtract \eqref{eq:DiscVarTINS} from \eqref{eq:VarTINS} and use arbitrary $\rb{\vv_h,q_h}\in\VV_h\times Q_h$ as test functions.
\end{thmProof}

\subsection{Conservation properties}	 \label{sec:ConservationProperties}
This section is based on \cite{CharnyiEtAl17} and directly attaches to Remark \ref{rem:ConsTrilinearForm}. It is well-known that most of the typical discretisations of the Navier--Stokes equations enforce the divergence-free condition only in a weak sense, in general leading to discrete velocities with $\nabla\cdot\uu_h\neq 0$; cf. \cite{JohnEtAl16}. As a consequence, even though the physics of the PDEs dictate otherwise, in addition to the potential lack of mass conservation, the following quantities also may not be conserved on the discrete level:
\begin{subequations}  \label{eq:ConsQuantities}
\begin{empheq}{alignat=2}
	\text{Kinetic energy:}\qquad 
		&E\rb{\uu,t} = \frac{1}{2}\int_\Omega \abs{\uu\rb{t,\x}}^2 \dx \qquad
			&&\text{for }\nu=0,~\ff=\zero	\label{eq:KinEn}\\
	\text{Linear momentum:}\qquad 
		&\MM\rb{\uu,t} = \int_\Omega \uu\rb{t,\x} \dx \qquad
			&&\text{for } \ff \text{ with zero linear momentum } \label{eq:LinMom}	\\
	\text{Angular momentum:}\qquad 
		&\MM_\x\rb{\uu,t} = \int_\Omega \uu\rb{t,\x}\times\x \dx	 \qquad
			&&\text{for } \ff \text{ with zero angular momentum } \label{eq:AngMom}
\end{empheq}
\end{subequations}
Let us, for a moment, drop the assumption of using a weakly divergence-free FEM. Then, in general, \eqref{eq:TrilinearSkewSymm} does not hold anymore and the particular choice of the trilinear form in the variational formulation is vital. Writing the momentum balance of the strong formulation with a generic nonlinear term $N\rb{\cdot}$ as
\begin{align}
	\partial_t\uu - \nu\Delta \uu + N\rb{\uu} +\nabla p = \ff,
\end{align}
we distinguish the following formulations:
\begin{subequations} \label{eq:NonlinearForms}
	\begin{align}
		\text{convective:} \qquad 
			&N\rb{\uu}=\rb{\uu\cdot\nabla}\uu \\
		\text{skew-symmetric:} \qquad 
			&N\rb{\uu}=\rb{\uu\cdot\nabla}\uu + \frac{1}{2}\rb{\nabla\cdot\uu}\uu \\
		\text{EMAC:} \qquad 
			&N\rb{\uu}=2\DD{\uu}\uu + \rb{\nabla\cdot\uu}\uu \label{eq:EMAC}
	\end{align}
\end{subequations}
For \eqref{eq:EMAC}, called `energy momentum and angular momentum conserving formulation', we use the decomposition of the velocity gradient tensor into deformation (symmetric) and spin (skew-symmetric) tensor:
\begin{align}
\nabla\uu	=\DD{\uu}+\SP{\uu} 
			= \frac{1}{2}\rb{\nabla\uu+\nabla\uu^\dag}+\frac{1}{2}\rb{\nabla\uu-\nabla\uu^\dag}	
\end{align}

\begin{thmThe}[Conservation properties \cite{CharnyiEtAl17}]
Regarding the corresponding FEM schemes to \eqref{eq:NonlinearForms}, the EMAC formulation conserves \eqref{eq:KinEn}-\eqref{eq:AngMom} whereas the skew-symmetric formulation only conserves \eqref{eq:KinEn} and the convective formulation does not conserve anything. However, for exactly divergence-free FEMs all nonlinear formulations are equivalent and therefore, all physical quantities \eqref{eq:ConsQuantities} of the discrete solution are conserved automatically.
\end{thmThe}

\section{Stability analysis and well-posedness}	\label{sec:StabWellPosed}
In this section, we derive stability estimates for the semi-discrete, continuous-in-time velocity approximation. Furthermore, we address the existence and uniqueness of the discrete variational formulation \eqref{eq:DiscVarTINS}. The analysis in this section is in some places inspired by \cite{BurmanFernandez07,ArndtEtAl15,DallmannEtAl16}.\\

Motivated by symmetric testing in the bi- and trilinear forms of \eqref{eq:DiscVarTINS}, that is
\begin{align}
\nu\norm{\nabla\vv}_\LP{2}^2=
a\rb{\vv,\vv}=
a\rb{\vv,\vv}+t\rb{\uu_h;\vv,\vv}+b\rb{\vv,q}-b\rb{\vv,q},
\end{align}
we define the following dissipation energy norm: 
\begin{equation} \label{eq:energynorm}
	\enorm{\vv}^2=\nu\norm{\nabla\vv}_\LP{2}^2=\nu\rb{\nabla\vv,\nabla\vv}
	=\nu\abs{\vv}_\HM{1}^2=a\rb{\vv,\vv}
\end{equation}

Our intent in the following is to decouple velocity and pressure. This can be done thanks to the LBB stability which implies $\VV_h^\dvg\neq\set{\zero}$. Indeed, from the second equation of \eqref{eq:DiscVarTINS} we know that  $\uu_h\in\VV_h^\dvg\subset\VV^\dvg$ for almost every $t\in\sqb{0,\tend}$, provided we additionally assume that $\uu_{0h}\in\VV_h^\dvg$. When we test the first equation of \eqref{eq:DiscVarTINS} only with discretely divergence-free functions, the discrete velocity solution can be computed from the following reduced, pressure-free variational formulation:
\begin{subequations}\label{eq:DiscinVdiv}
	\begin{empheq}[left=\empheqlbrace]{align} 
	&\text{find }\uu_h\colon\rb{0,\tend}\to\VV_h^\dvg
	\text{ with }\uu_h\rb{0}=\uu_{0h} \text{ s.t. }\forall\vv_h\in\VV_h^\dvg\\
		&\rb{\partial_t\uu_h,\vv_h} + a\rb{\uu_h,\vv_h} + t\rb{\uu_h;\uu_h,\vv_h} =\bra{\ff,\vv_h}		
	\end{empheq} 
\end{subequations}

The following result shows that for the discrete solution, the kinetic energy at every time instance $t\in\sqb{0,\tend}$ and the total kinetic energy dissipated over $\sqb{0,t}$ is bounded by the data. 

\begin{thmLem}[Energy estimate] \label{lem:EnEst}
Let $\ff\in\Lp{2}\rb{0,\tend;\VV^*}$ and $\uu_{0h}\in\VV_h^\dvg$. Then, for each $0\leqslant t\leqslant \tend$ we obtain
\begin{equation}
\norm{\uu_h\rb{t}}_\LP{2}^2 + \int_0^t \enorm{\uu_h\rb{\tau}}^2\dtau \leqslant
\norm{\uu_{0h}}_\LP{2}^2 + \frac{1}{\nu}\norm{\ff}_{\Lp{2}\rb{0,\tend;\VV^*}}^2.
\end{equation}
\end{thmLem}

\begin{thmProof}
Cf., for example, \cite[Section 9.2]{Layton08}.	
\end{thmProof}

\begin{thmRem}
If we assume the more restrictive regularity $\ff\in\Lp{1}\rb{0,\tend;\LP{2}}$	in Lemma \ref{lem:EnEst} one can even remove the explicit dependence on the viscosity and prove \cite[Lemma 3.1]{ArndtEtAl15}
\begin{equation} \label{eq:EnergyMoreRegf}
\frac{1}{2}\norm{\uu_h\rb{t}}_\LP{2}^2 + \int_0^t \enorm{\uu_h\rb{\tau}}^2\dtau \leqslant
\norm{\uu_{0h}}_\LP{2}^2 + \frac{3}{2}\norm{\ff}_{\Lp{1}\rb{0,\tend;\LP{2}}}^2.
\end{equation} 
\end{thmRem}

\begin{thmCor}[Existence and uniqueness]
If $\ff$ is Lipschitz continuous in time, there exists a unique semi-discrete velocity $\uu_h\colon\sqb{0,\tend}\to\VV_h^\dvg$ of the variational formulation \eqref{eq:DiscinVdiv}.	
\end{thmCor}

\begin{thmProof}
Cf., for example, \cite[Section 9.2]{Layton08}. On the right-hand side, the semi-discrete problem	
\begin{equation} \label{eq:SysODE}
	\rb{\partial_t \uu_h,\vv_h} =  \bra{\ff,\vv_h}-a\rb{\uu_h,\vv_h} - t\rb{\uu_h;\uu_h,\vv_h} 
\end{equation}
has a quadratic and thus locally Lipschitz continuous nonlinearity. As a consequence of Lemma \ref{lem:EnEst}, for a fixed $\nu>0$, every potential solution of \eqref{eq:SysODE} cannot blow up in finite time and, therefore, there exists a unique and global-in-time solution due to the theory of first-order ODEs.
\end{thmProof}

\section{Pressure- and semi-robust error estimates}	\label{sec:ErrorEstimatesVel}

In this section, we assume that $\rb{\uu,p}\in\VV\times\Q$ solves \eqref{eq:VarTINS} and $\rb{\uu_h,p_h}\in\VV_h\times\Q_h$ is the unique FE solution to \eqref{eq:DiscVarTINS}. At some places, ideas from \cite{BurmanFernandez07,ArndtEtAl15,DallmannEtAl16} are used again.  With appropriate approximation operators $\rb{\PI_h,\pi_h}\colon\VV\times\Q\to\VV_h\times\Q_h$ we decompose the error as
\begin{subequations}
\begin{align}
\uu-\uu_h&=\rb{\uu-\PI_h\uu}+\rb{\PI_h\uu-\uu_h}=\uETA + \uERR =\uXI	\\
p-p_h&=\rb{p-\pi_h p}+\rb{\pi_h p -p_h}=\pETA+\pERR=\pXI
\end{align}	
\end{subequations}
and refer to $\rb{\uETA,\pETA}$ and $\rb{\uERR,\pERR}$ as approximation and discretisation error, respectively.

\subsection{Discrete Stokes projection}	
For the above error splitting, we introduce here the discrete Stokes projection which is used frequently for  numerical analysis \cite{ErnGuermond04,Ingram13}. The stationary incompressible Stokes problem with no-slip boundary conditions and general body force term $\gbld$ is given by
\begin{subequations}\label{eq:Stokes}
	\begin{empheq}[left=\empheqlbrace]{alignat=2} 
		-\nu\Delta \uu_s +\nabla p_s &=  \gbld \qquad 		&&\text{in }\Omega, 			\\
		\nabla\cdot\uu_s &= 0 								&&\text{in }\Omega, 			\\
		\uu_s &= \zero 										&&\text{on }\partial\Omega,	
	\end{empheq} 
\end{subequations}
for which the standard finite element formulation reads as follows:
\begin{subequations}\label{eq:VarStokes}
	\begin{empheq}[left=\empheqlbrace]{align} 
	\text{find }\rb{\uu_{sh},p_{sh}}\in\VV_h\times\Q_h\text{ s.t. }&\forall\rb{\vv_h,q_h}\in\VV_h\times\Q_h	\\
		 \nu\rb{\nabla\uu_{sh},\nabla\vv_h}  - \rb{p_{sh},\nabla\cdot\vv_h}  &= \rb{\gbld,\vv_h}  \\
		 \rb{q_h,\nabla\cdot\uu_{sh}} & = 0
	\end{empheq} 
\end{subequations}
For this Galerkin approximation, error estimates are well-known \cite{GiraultRaviart86,ErnGuermond04}:
\begin{subequations} \label{eq:ErrStokes}
\begin{align}
	\norm{\uu_s-\uu_{sh}}_\LP{2} 
		+ h\abs{\uu_s-\uu_{sh}}_\HM{1} 
		&\lesssim h \rb{\inf_{\vv_h\in\VV_h}\norm{\uu_s-\vv_h}_\HM{1}
		+  \nu^{-1}\inf_{q_h\in\Q_h}\norm{p_s-q_h}_\Lp{2}} \\
	\norm{p_s-p_{sh}}_\Lp{2}
		&\lesssim \nu \inf_{\vv_h\in\VV_h}\norm{\uu_s-\vv_h}_\HM{1}
		+  \inf_{q_h\in\Q_h}\norm{p_s-q_h}_\Lp{2}	
\end{align}	
\end{subequations}

\begin{thmRem}
In order to obtain the estimate for $\norm{\uu_s-\uu_{sh}}_\LP{2}$ an Aubin--Nitsche duality argument is mandatory. However, to be able to apply such an argument, it is necessary for the Stokes problem to have smoothing properties. Referring to \cite[Section 4.2]{ErnGuermond04}, a sufficient but rather restrictive assumption guaranteeing this is that $\Omega$ is either a convex polygon in $d=2$ or $\Omega$ is of class $\mathcal{C}^{1,1}$ for $d\in\set{2,3}$.
\end{thmRem}
Following the concept of \cite{FrutosEtAl16}, the forcing term is chosen in particular as
\begin{align}
	\gbld = \ff -\partial_t\uu - \rb{\uu\cdot\nabla}\uu -\nabla p = -\nu\Delta\uu,
\end{align}
where $\rb{\uu,p}\in\VV\times\Q$ denotes the exact solution of the Navier--Stokes problem \eqref{eq:TINS}. Therefore, on the continuous level, we are searching for the unique solution $\rb{\uu_s,p_s}$ fulfilling 
\begin{equation}
	-\nu\Delta \uu_s +\nabla p_s =  \ff -\partial_t\uu - \rb{\uu\cdot\nabla}\uu -\nabla p.
\end{equation}
Since $\rb{\uu,p}$ solves the Navier--Stokes equations uniquely, the exact solution of the Stokes equations \eqref{eq:Stokes} with $\gbld=-\nu\Delta\uu$ is given by $\rb{\uu_s,p_s}=\rb{\uu,0}$. This leads to the following definition.

\begin{thmDef}[Discrete Stokes projection] \label{def:StokesProjection}
Let $\rb{\uu,p}\in\VV\times\Q$ solve \eqref{eq:VarTINS}. Then, the discrete Stokes projection $\rb{\PIs \uu,\pis p}\in\VV_h\times\Q_h$ is defined as the unique solution $\rb{\uu_{sh},p_{sh}}$ to \eqref{eq:VarStokes} with $\gbld=-\nu\Delta\uu$.
\end{thmDef}

 In this way, the approximation properties of the projection operators can be derived from the error estimates for the Stokes problem \eqref{eq:ErrStokes}:
\begin{subequations} \label{eq:ErrStokesProj}
\begin{align}
	\norm{\uu-\PIs\uu}_\LP{2} 
		+ h\abs{\uu-\PIs\uu}_\HM{1} 
		&\leqslant C_s h \inf_{\vv_h\in\VV_h}\norm{\uu-\vv_h}_\HM{1} \label{eq:ErrStokesProjVel} \\
	\norm{\pis p}_\Lp{2}
		&\leqslant C_s \nu \inf_{\vv_h\in\VV_h}\norm{\uu-\vv_h}_\HM{1}		
\end{align}	
\end{subequations}

Furthermore, by construction, the Stokes projection is discretely divergence-free, that is $\PIs\uu\in\VV_h^\dvg$. Assuming sufficient smoothness in time such that $\gbld=\partial_t\rb{-\nu\Delta\uu}$ can be chosen in \eqref{eq:VarStokes}, one obtains 
\begin{equation} \label{eq:ErrDtStokesProj}
\norm{\partial_t\rb{\uu-\PIs\uu}}_\LP{2} 
		+ h\abs{\partial_t\rb{\uu-\PIs\uu}}_\HM{1} 
		\leqslant C_s h \inf_{\vv_h\in\VV_h}\norm{\partial_t\uu-\vv_h}_\HM{1}	.
\end{equation}

\begin{thmCor} \label{cor:StokesProj}
Let $\rb{\uu,p}\in\VV\times\Q$ solve \eqref{eq:VarTINS} and $\rb{\PI_h\uu,\pi_h p}=\rb{\PIs\uu,\pis p}\in\VV_h^\dvg\times\Q_h$ be the Stokes projection of Definition \ref{def:StokesProjection}. Then, the following holds:
\begin{enumerate}[label=(\roman*)]
\item Error equation:
	\begin{equation}
		a\rb{\uETA,\vv_h} - b\rb{\vv_h, \pis p} -b\rb{\uETA,q_h} = 0,\quad \forall\rb{\vv_h,q_h}\in\VV_h\times\Q_h
	\end{equation}
\item Stokes and Ritz projection coincide on $\VV_h^\dvg$, that is 
	\begin{equation}
	\rb{\nabla \PIs\uu,\nabla\vv_h} = \rb{\nabla \uu,\nabla\vv_h},\quad\forall\vv_h\in\VV_h^\dvg.	
	\end{equation}
\item Max-norm stability: Assuming additionally that $\Omega$ is convex, we have
	\begin{equation} \label{eq:StokesProjMaxStab}
		\abs{\PIs\uu}_\WMP{1}{\infty}
			\leqslant C_\infty \abs{\uu}_\WMP{1}{\infty}.
	\end{equation}
\end{enumerate}
\end{thmCor}

\begin{thmProof}
Ad (i): Since $\uu\in\VV$ solves \eqref{eq:VarTINS}, we know that $\uu\in\VV^\dvg$ and thus, $b\rb{\uu,q_h}=0$  for all $q_h\in\Q_h$. Thus, inserting $\gbld=-\nu\Delta\uu$ in \eqref{eq:VarStokes} and integrating by parts yields the assertion:
 \begin{subequations}
 	\begin{align}
 		\nu\rb{\nabla\PIs\uu,\nabla\vv_h}  - \rb{\pis p,\nabla\cdot\vv_h}  &= \nu\rb{\nabla\uu,\nabla\vv_h} \\
		 \rb{q_h,\nabla\cdot\PIs\uu} & = \rb{q_h,\nabla\cdot\uu}
 	\end{align}
 \end{subequations}
 
Ad (ii):  Now, suppose that we test in (i) with $\vv_h\in\VV_h^\dvg$. Then, by construction, $- b\rb{\vv_h, \pis p}=0$. 

Ad (iii): In general, for an approach from the direction of the Ritz projection, see \cite{GuzmanEtAl09}.  In order to obtain the result directly from the Stokes equations, we refer to \cite{GiraultEtAl15}; notice, however, that in the present work we define the Stokes projection in a slightly different way. Having said this, it is nevertheless possible to modify the approach in \cite{GiraultEtAl15} to obtain a pressure-independent max-norm estimate whenever a weakly divergence-free FEM is used to compute the Stokes projection.
\end{thmProof}

\begin{thmRem} \label{rem:b(v_h,p)=0}
Note that in this subsection we did not use the inclusion property \eqref{eq:divVhinQh} and, therefore, all previous explanations also hold true verbatim for general inf-sup stable and conforming FEMs. However, in order to separate the pressure from the velocity error estimates, we will need the property $\nabla\cdot\VV_h\subset\Q_h$, or more precisely $\VV_h^\dvg\subset\VV^\dvg$, for weakly divergence-free FEM to infer that $b\rb{\vv_h, p}=0$ for all $\vv_h\in\VV_h^\dvg$.
\end{thmRem}

\begin{thmRem}
The choice of an appropriate approximation operator for the error splitting is a crucial part of the numerical analysis. Note that in \cite{ArndtEtAl15,DallmannEtAl16} the error analysis is carried out in $\VV_h^\dvg$ with the use of a divergence-preserving interpolation operator \cite{GiraultScott03}. 
\end{thmRem}

\subsection{Unstabilised Galerkin-FEM}	

Now, let $\rb{\PI_h\uu,\pi_h p}=\rb{\PIs\uu,\pis p}\in\VV_h^\dvg\times\Q_h$ be the Stokes projection of Definition \ref{def:StokesProjection}.

\begin{thmLem}[Difference of convective terms] \label{lem:ConvTerms}
Let $\uu\in\VV$ solve \eqref{eq:VarTINS}, $\uu_h\in\VV_h$ solve \eqref{eq:DiscVarTINS} and suppose that $\uu\in\Lp{\infty}\rb{0,\tend;\WMP{1}{\infty}}$. Then, for all finite $\eps>0$ and $\sigma>0$, the following estimate holds true:
\begin{subequations}
\begin{align}
|t\rb{\uu;\uu,\uERR}&-t\rb{\uu_h;\uu_h,\uERR}| \leqslant 
\frac{1}{4\eps\sigma^2}\norm{\uETA}_\LP{2}^2
+ \frac{1}{4\eps}\abs{\uETA}_\HM{1}^2  \\ &+
\rb{2\abs{\uu}_\WMP{1}{\infty}
+\eps \sigma^2\abs{\uu}_\WMP{1}{\infty}^2
+\eps \norm{\PIs\uu}_\LP{\infty}^2
+ \abs{\PIs\uu}_\WMP{1}{\infty}}\norm{\uERR}_\LP{2}^2 
\end{align}	
\end{subequations}
\end{thmLem}

\begin{thmProof}
Using $\uu_h=\PIs\uu-\uERR$, we deduce that
\begin{align}
t\rb{\uu_h;\uu_h,\uERR}=\rb{\uu_h\cdot\nabla\uu_h,\uERR}	 
= \rb{\uu_h\cdot\nabla \PIs\uu,\uERR} - \rb{\uu_h\cdot\nabla \uERR,\uERR}
=\rb{\uu_h\cdot\nabla \PIs\uu,\uERR}, 
\end{align}
where the last equality is due to Lemma \ref{lem:Trilinear} together with $\nabla\cdot\uu_h=0$. Therefore, we obtain
\begin{subequations}
\begin{align}
t\rb{\uu;\uu,\uERR}-t\rb{\uu_h;\uu_h,\uERR} &= 
t\rb{\uu;\uu,\uERR}-t\rb{\uu_h;\uu,\uERR}+t\rb{\uu_h;\uu,\uERR}-\rb{\uu_h\cdot\nabla \PIs\uu,\uERR} 	 \\
&=\int_\Omega\rb{\sqb{\uu-\uu_h}\cdot\nabla}\uu\cdot\uERR\dx + \int_\Omega\rb{\uu_h\cdot\nabla} \sqb{\uu-\PIs\uu}\cdot \uERR\dx \\ 
&=t\rb{\uETA+\uERR;\uu,\uERR} + t\rb{\uu_h;\uETA,\uERR} =\fT_1 + \fT_2.
\end{align}	
\end{subequations}

Mainly using $\uu\rb{t}\in\WMP{1}{\infty}{\rb{\Omega}}$ for $0\leqslant t\leqslant T$ and Young's inequality, we can estimate these two terms separately. For the first one, with the estimate from Lemma \ref{lem:Trilinear} and arbitrary $\eps>0$ and $\sigma> 0$, we get 
\begin{subequations}
\begin{align}
\abs{\fT_1}&=\abs{t\rb{\uETA+\uERR;\uu,\uERR}} \leqslant
\rb{\norm{\uETA}_\LP{2}+\norm{\uERR}_\LP{2}}\norm{\nabla\uu}_\LP{
\infty}\norm{\uERR}_\LP{2} \\ &=
\abs{\uu}_\WMP{1}{\infty}\norm{\uERR}_\LP{2}^2 
+ \sigma^{-1}\norm{\uETA}_\LP{2}\sigma\abs{\uu}_\WMP{1}{\infty}\norm{\uERR}_\LP{2} \\ &\leqslant
\frac{1}{4\eps \sigma^2}\norm{\uETA}_\LP{2}^2
+\rb{\abs{\uu}_\WMP{1}{\infty}+\eps \sigma^2\abs{\uu}_\WMP{1}{\infty}^2}\norm{\uERR}_\LP{2}^2. 
\end{align}	
\end{subequations}

For the second term, again using $\uu_h=\PIs\uu-\uERR$ and the triangle inequality, we deduce that
\begin{align}
	\abs{\fT_2}=\abs{t\rb{\uu_h;\uETA,\uERR}}&\leqslant
	\abs{t\rb{\PIs\uu;\uETA,\uERR}} + \abs{t\rb{\uERR;\uETA,\uERR}}
	=\abs{\fT_{2,1}}+\abs{\fT_{2,2}},
\end{align}

and, for the first part of the partition, with Young's inequality ($\eps>0$) obtain
\begin{align}
	\abs{\fT_{2,1}}&=\abs{t\rb{\PIs\uu;\uETA,\uERR}} \leqslant
	\norm{\PIs\uu}_\LP{\infty}\norm{\nabla\uETA}_\LP{2}\norm{\uERR}_\LP{2} 
	\leqslant \eps \norm{\PIs\uu}_\LP{\infty}^2\norm{\uERR}_\LP{2}^2
		+ \frac{1}{4\eps}\abs{\uETA}_\HM{1}^2
\end{align}	

and, by virtue of Lemma \ref{lem:Trilinear} plus the triangle inequality,
\begin{subequations}
\begin{align}
	\abs{\fT_{2,2}}&=\abs{t\rb{\uERR;\uETA,\uERR}} \leqslant 
	\norm{\nabla\uETA}_\LP{\infty} \norm{\uERR}_\LP{2}^2 
	= \abs{\uETA}_\WMP{1}{\infty} \norm{\uERR}_\LP{2}^2 \\
	&\leqslant \rb{\abs{\uu}_\WMP{1}{\infty}+\abs{\PIs\uu}_\WMP{1}{\infty}} 
		\norm{\uERR}_\LP{2}^2. 
\end{align}	
\end{subequations}

Combining the three estimates concludes the proof.
\end{thmProof}

\begin{thmThe}[Semi-robust velocity discretisation error estimate] \label{thm:VelDiscErr}
Let $\uu\in\VV$ solve \eqref{eq:VarTINS}, $\uu_h\in\VV_h$ solve \eqref{eq:DiscVarTINS} and assume that $\uu\in\Lp{\infty}\rb{0,\tend;\WMP{1}{\infty}}$, $\partial_t\uu\in\Lp{2}{\rb{0,\tend;\LP{2}}}$ and $\uu_h\rb{0}=\PIs \uu_0$. Then, we obtain the following error estimate for the velocity of the FEM:
\begin{align} \label{eq:uErrEst}
\norm{\uERR}_{\Lp{\infty}\rb{0,\tend;\LP{2}}}^2 &+ \int_0^\tend \enorm{\uERR\rb{\tau}}^2 \dtau \leqslant 
e^{C_\uu\tend}\int_0^\tend \sqb{\norm{\partial_t\uETA\rb{\tau}}_\LP{2}^2+\frac{1}{h}\norm{\uETA\rb{\tau}}_\LP{2}^2+ \abs{\uETA\rb{\tau}}_\HM{1}^2} \dtau
\end{align}		
Here, the Gronwall constant is given by
\begin{align}
	C_\uu=
1&+4\abs{\uu}_{\Lp{\infty}\rb{0,\tend;\WMP{1}{\infty}}}
+ h\abs{\uu}_{\Lp{\infty}\rb{0,\tend;\WMP{1}{\infty}}}^2
+ \norm{\PIs\uu}_{\Lp{\infty}\rb{0,\tend;\LP{\infty}}}^2
+ 2\abs{\PIs\uu}_{\Lp{\infty}\rb{0,\tend;\WMP{1}{\infty}}}.
\end{align}	
\end{thmThe}

\begin{thmRem}
	Assuming additionally that $\uu\in\Lp{2}{\rb{0,\tend;\HM{1}}}$, Theorem \ref{thm:VelDiscErr} directly implies strong velocity convergence of the FEM; cf. \cite{ArndtEtAl15}.
\end{thmRem}

\begin{thmProof}
Use $\rb{\uERR,\pERR}\in\VV_h^\dvg\times\Q_h$ as test functions in Corollary \ref{cor:GalOrtho}: 
\begin{align}
0&=\rb{\partial_t\uXI,\uERR}+a\rb{\uXI,\uERR}+t\rb{\uu;\uu,\uERR}-t\rb{\uu_h;\uu_h,\uERR}+b\rb{\uERR,\pXI} -b\rb{\uXI,\pERR}
\end{align}

Inserting the definition of the full errors yields
\begin{subequations}
\begin{align}
	0=\rb{\partial_t\uERR,\uERR}&+a\rb{\uERR,\uERR}+a\rb{\uETA,\uERR}
	+b\rb{\uERR,p}-b\rb{\uERR,\pis p}+b\rb{\uERR,\pERR} -b\rb{\uERR,\pERR}\\ 
	&-b\rb{\uETA,\pERR}+\rb{\partial_t\uETA,\uERR} +\sqb{t\rb{\uu;\uu,\uERR} - t\rb{\uu_h;\uu_h,\uERR}}.
\end{align}	
\end{subequations}

Now, using the properties of the Stokes projection, namely Corollary \ref{cor:StokesProj} and Remark \ref{rem:b(v_h,p)=0}, a lot of terms drop out and, after rearranging, we obtain
\begin{align}
	\rb{\partial_t\uERR,\uERR}+a\rb{\uERR,\uERR} =
	-\rb{\partial_t\uETA,\uERR}-\sqb{t\rb{\uu;\uu,\uERR} - t\rb{\uu_h;\uu_h,\uERR}}.
\end{align}	

Note that we specifically used the property $\VV_h^\dvg\subset\VV^\dvg$ for weakly divergence-free FEM to infer that $b\rb{\uERR, p}=0$ and, therefore, the pressure drops out of the velocity error estimates entirely. Bearing in mind that $\rb{\partial_t\uERR,\uERR}=\frac{1}{2}\frac{\drm}{\drm t}\norm{\uERR}_\LP{2}^2$ and using the definition of the energy norm \eqref{eq:energynorm} leads to
\begin{align}
	\frac{1}{2}\frac{\drm}{\drm t}\norm{\uERR}_\LP{2}^2 + \enorm{\uERR}^2 = 
	-\rb{\partial_t\uETA,\uERR} -\sqb{t\rb{\uu;\uu,\uERR} - t\rb{\uu_h;\uu_h,\uERR}}.
\end{align}

The convective terms have been estimated in Lemma \ref{lem:ConvTerms} for all finite $\eps>0$ and $\sigma>0$:
\begin{subequations}
\begin{align}
|t\rb{\uu;\uu,\uERR}&-t\rb{\uu_h;\uu_h,\uERR}| \leqslant 
\frac{1}{4\eps\sigma^2}\norm{\uETA}_\LP{2}^2
+ \frac{1}{4\eps}\abs{\uETA}_\HM{1}^2  \\ &+
\rb{2\abs{\uu}_\WMP{1}{\infty}
+\eps \sigma^2\abs{\uu}_\WMP{1}{\infty}^2
+\eps \norm{\PIs\uu}_\LP{\infty}^2
+ \abs{\PIs\uu}_\WMP{1}{\infty}}\norm{\uERR}_\LP{2}^2 
\end{align}	
\end{subequations}

Assuming that $\partial_t\uETA \in\LP{2}{\rb{\Omega}}$, the dynamical first term on the right-hand side can be estimated with Cauchy--Schwarz and Young's inequality:
\begin{align}
\abs{\rb{\partial_t\uETA,\uERR}}\leqslant\norm{\partial_t\uETA}_\LP{2}\norm{\uERR}_\LP{2} \leqslant \frac{1}{2}\norm{\partial_t\uETA}_\LP{2}^2	+\frac{1}{2}\norm{\uERR}_\LP{2}^2	
\end{align}

Combining all above estimates yields
\begin{subequations}
\begin{align}
	\frac{1}{2}\frac{\drm}{\drm t}&\norm{\uERR}_\LP{2}^2 
	+ \enorm{\uERR}^2 \leqslant
		\frac{1}{2}\norm{\partial_t\uETA}_\LP{2}^2 +
	\frac{1}{4\eps\sigma^2}\norm{\uETA}_\LP{2}^2
+ \frac{1}{4\eps}\abs{\uETA}_\HM{1}^2
	 \\ &+ \rb{\frac{1}{2}+2\abs{\uu}_\WMP{1}{\infty}
+\eps \sigma^2\abs{\uu}_\WMP{1}{\infty}^2
+\eps \norm{\PIs\uu}_\LP{\infty}^2
+ \abs{\PIs\uu}_\WMP{1}{\infty}}\norm{\uERR}_\LP{2}^2.
\end{align}	
\end{subequations}

Choosing $\sigma=\sqrt{h}$, $\eps=\nicefrac{1}{2}$ and multiplying this inequality by 2 results in 
\begin{subequations}\label{eq:Gronwall1}
\begin{align}
	\frac{\drm}{\drm t}&\norm{\uERR}_\LP{2}^2 
	+ 2\enorm{\uERR}^2 \leqslant
		\norm{\partial_t\uETA}_\LP{2}^2 +
	\frac{1}{h}\norm{\uETA}_\LP{2}^2+\abs{\uETA}_\HM{1}^2
	 \\ &+ \rb{1+4\abs{\uu}_\WMP{1}{\infty}
+ h\abs{\uu}_\WMP{1}{\infty}^2
+ \norm{\PIs\uu}_\LP{\infty}^2
+ 2\abs{\PIs\uu}_\WMP{1}{\infty}}\norm{\uERR}_\LP{2}^2.
\end{align}	
\end{subequations}

Using the abbreviation
\begin{equation}
g\rb{\tau} = 1+4\abs{\uu\rb{\tau}}_\WMP{1}{\infty}
+ h\abs{\uu\rb{\tau}}_\WMP{1}{\infty}^2
+ \norm{\PIs\uu\rb{\tau}}_\LP{\infty}^2
+ 2\abs{\PIs\uu\rb{\tau}}_\WMP{1}{\infty},	
\end{equation}

from the regularity assumption, we infer that for all $0\leqslant t\leqslant T$
\begin{equation}
	G\rb{t}=\int_0^t g\rb{\tau}\dtau < \infty.
\end{equation}

By applying Gronwall's lemma \cite[Lemma 6.9]{ErnGuermond04} and acknowledging the monotonicity of $G$, we can drop the exponential term on the left-hand side and, therefore, for all $0\leqslant t \leqslant T$, reduce our inequality to 
\begin{align}
\norm{\uERR\rb{t}}_\LP{2}^2 &+ \int_0^t \enorm{\uERR\rb{\tau}}^2 \dtau \leqslant
\int_0^t e^{G\rb{t}-G\rb{\tau}}\sqb{\norm{\partial_t\uETA\rb{\tau}}_\LP{2}^2+\frac{1}{h}\norm{\uETA\rb{\tau}}_\LP{2}^2
+ \abs{\uETA\rb{\tau}}_\HM{1}^2} \dtau.
\end{align}	

On the other hand, using the regularity assumption $\uu\in\Lp{\infty}\rb{0,\tend;\WMP{1}{\infty}}$, we can estimate
\begin{equation}
	G\rb{t}-G\rb{\tau}=
	\int_\tau^t g\rb{s} \drm s
	\leqslant C_\uu\sqb{t-\tau}
	\leqslant C_\uu\tend
\end{equation}

with a Gronwall constant
\begin{align}
	C_\uu=
1&+4\abs{\uu}_{\Lp{\infty}\rb{0,\tend;\WMP{1}{\infty}}}
+ h\abs{\uu}_{\Lp{\infty}\rb{0,\tend;\WMP{1}{\infty}}}^2 
+ \norm{\PIs\uu}_{\Lp{\infty}\rb{0,\tend;\LP{\infty}}}^2
+ 2\abs{\PIs\uu}_{\Lp{\infty}\rb{0,\tend;\WMP{1}{\infty}}}.
\end{align}	

Given \eqref{eq:StokesProjMaxStab}, this concludes the proof.
\end{thmProof}

\begin{thmRem}
Note that the error estimate \eqref{eq:uErrEst} for the velocity is independent of the pressure, that is, pressure-robust. This is characteristic for weakly divergence-free methods and a major advantage over other FEM. For a more general discussion of pressure-robustness, we refer to \cite{LinkeEtAl16,JohnEtAl16}.	
\end{thmRem}

\begin{thmRem}
It is possible to relax the regularity assumptions in Theorem \ref{thm:VelDiscErr}. In fact, we only require that $G\rb{t}<\infty$ for all $0\leqslant t \leqslant \tend$. Thus, assuming $\uu,\PIs\uu\in\Lp{2}{\rb{0,\tend;\WMP{1}{\infty}}}$ would be sufficient. However, the given regularity seems to be rather natural to us. Furthermore, we would like to make reference to \cite{Burman15}, where the analysis for the two-dimensional Navier--Stokes problem is based on a scale separation $\uu=\overline{\uu}+\uu^\prime$ into large eddies $\overline{\uu}$ and small scales (fluctuations) $\uu^\prime$. Assuming that $\overline{\uu}\in\WMP{1}{\infty}$ and $\uu^\prime\in\WMP{1}{p}\cap\LP{\infty}$ for $p>d=2$, this leads to error estimates where the Gronwall constant $\exp\rb{1+C\norm{\nabla\uu}_{\Lp{\infty}\rb{0,\tend;\LP{\infty}}}+\dots}$ can be replaced by $\exp\rb{1+C\norm{\nabla\overline{\uu}}_{\Lp{\infty}\rb{0,\tend;\LP{\infty}}}+\dots}$, which is much more realistic.

\end{thmRem}

\begin{thmCor}[Velocity discretisation error convergence rate] \label{cor:VelConvRate}
Under the assumptions of the previous theorem, assume a smooth solution according to 
\begin{align}
	\uu&\in\Lp{\infty}\rb{0,\tend;\WMP{1}{\infty}}\cap\Lp{2}\rb{0,\tend;\HM{r}}, \quad  
	\partial_t\uu\in\Lp{2}\rb{0,\tend;\HM{r-1}}.	
\end{align}
Then, we obtain the following convergence rate for the velocity of the FEM:
\begin{align}
\norm{\uERR}_{\Lp{\infty}\rb{0,\tend;\LP{2}}}^2 &+ \int_0^\tend \enorm{\uERR\rb{\tau}}^2 \dtau 
	\leqslant Ch^{2\rb{r_\uu-1}}e^{C_\uu\tend} 
		\int_0^\tend \sqb{\rb{1+h}\abs{\uu\rb{\tau}}_\HM{r_\uu}^2
		+\abs{\partial_t\uu\rb{\tau}}_\HM{r_\uu-1}^2} \dtau \label{eq:VelDiscConvRate}
\end{align}		
Here, $r_\uu=\min\set{r,k+1}$ and the Gronwall constant $C_\uu$ is given in Theorem \ref{thm:VelDiscErr}. 
\end{thmCor}

\begin{thmProof}
The proof is a consequence of the approximation result \eqref{eq:ErrStokesProjVel} and \eqref{eq:ErrDtStokesProj} for the Stokes projection and the interpolation result \eqref{eq:uInt} for the velocity space:
\begin{subequations}
\begin{align}
	\frac{1}{h}\norm{\uETA\rb{\tau}}_\LP{2}^2	
		&\leqslant C_s^2 h \inf_{\vv_h\in\VV_h}\norm{\uu\rb{\tau}-\vv_h}_\HM{1}^2
			\leqslant C h^{2r_\uu-1}\abs{\uu\rb{\tau}}_\HM{r_\uu}^2 \\
	\abs{\uETA\rb{\tau}}_\HM{1}^2 
		&\leqslant C_s^2 \inf_{\vv_h\in\VV_h}\norm{\uu\rb{\tau}-\vv_h}_\HM{1}^2
			\leqslant C h^{2\rb{r_\uu-1}}\abs{\uu\rb{\tau}}_\HM{r_\uu}^2 \label{eq:IntResultConv}\\
	\norm{\partial_t\uETA\rb{\tau}}_\LP{2}^2  
		&\leqslant C_s^2 h^2 \inf_{\vv_h\in\VV_h}\norm{\partial_t\uu\rb{\tau}-\vv_h}_\HM{1}^2
			\leqslant C h^{2\rb{r_\uu-1}}\abs{\partial_t\uu\rb{\tau}}_\HM{r_\uu-1}^2 \label{eq:IntResultDtu}
\end{align}	
\end{subequations}

Those are the particular terms in Theorem \ref{thm:VelDiscErr} which must be estimated.
\end{thmProof}

Note that, using the triangle inequality, \eqref{eq:uErrEst} and \eqref{eq:VelDiscConvRate} also hold true for the full error $\uXI=\uu-\uu_h$.

\begin{thmLem}[Max-norm estimate] \label{lem:uhInfinity}
Under the assumptions of the previous corollary, additionally assume $\uu\in\Lp{\infty}\rb{0,\tend;\HM{r}}$. Then,
\begin{subequations}
\begin{align}
	\norm{\uu-\uu_h}_{\Lp{\infty}\rb{0,\tend;\LP{\infty}}}^2
		&\leqslant Ch^2 \norm{\uu}_{\Lp{\infty}\rb{0,\tend;\WMP{1}{\infty}}}^2
		\\&+C h^{2\rb{r_\uu-1}-d} e^{C_\uu\tend} 
			\max\set{\abs{\uu}_{\Lp{2}\rb{0,\tend;\HM{r_\uu}}}^2,
			\abs{\partial_t\uu}_{\Lp{2}\rb{0,\tend;\HM{r_\uu-1}}}^2},
\end{align}	
\end{subequations}
where $r_\uu=\min\set{r,k+1}$. Especially, provided $r_\uu>\frac{d}{2}+1$, this shows that  $\uu_h\in\Lp{\infty}\rb{0,\tend;\LP{\infty}}$.
\end{thmLem}

\begin{thmProof}
	First of all, we introduce the orthogonal $\LTWO$-projection $\PI_0$ onto $\VV_h$. Then, as argued in \cite{BurmanFernandez07}, quasi-uniformity of the mesh is sufficient to obtain the following estimate:
	\begin{equation} \label{eq:PI0infinity}
		\norm{\uu-\PI_0\uu}_\LP{\infty}
			+h\norm{\uu-\PI_0\uu}_\WMP{1}{\infty}
			\leqslant Ch \norm{\uu}_\WMP{1}{\infty}
	\end{equation}
	
Now, using the triangle inequality, we conclude that
\begin{equation}
	\norm{\uu-\uu_h}_\LP{\infty}
		\leqslant \norm{\uu-\PI_0\uu}_\LP{\infty}
		+ \norm{\PI_0\uu-\uu_h}_\LP{\infty},
\end{equation}	
	
and the first term can be bounded by $Ch \norm{\uu}_\WMP{1}{\infty}$, see \eqref{eq:PI0infinity}. For the second term, we use a global inverse inequality (here, quasi-uniformity of the mesh is also required) based on the inverse inequality \eqref{eq:LocInvEq}  with $\ell=m=0$, $p=\infty$ and $q=2$, which yields
\begin{equation}
	\norm{\PI_0\uu-\uu_h}_\LP{\infty}
		\leqslant C_\inv h^{-\nicefrac{d}{2}} \norm{\PI_0\uu-\uu_h}_\LP{2}.	
\end{equation}

Once again with the triangle inequality, we deduce
\begin{equation} \label{eq:ProofMaxNorm}
	\norm{\PI_0\uu-\uu_h}_\LP{2}
		\leqslant \norm{\PI_0\uu-\uu}_\LP{2} + \norm{\uu-\uu_h}_\LP{2}.
\end{equation}

Using the approximation properties of smooth functions, see \cite[Proposition 1.134]{ErnGuermond04}, we estimate
\begin{align}
\norm{\PI_0\uu-\uu}_\LP{2}\leqslant Ch^{r_\uu}\abs{\uu}_\HM{r_\uu}	.
\end{align}

For the second term in \eqref{eq:ProofMaxNorm}, we invoke Corollary \ref{cor:VelConvRate}, which yields
\begin{align}
	\norm{\rb{\uu-\uu_h}\rb{t}}_\LP{2}
		\leqslant C h^{r_\uu-1} e^{C_\uu\nicefrac{\tend}{2}} \max\set{\abs{\uu}_{\Lp{2}\rb{0,\tend;\HM{r_\uu}}},\abs{\partial_t\uu}_{\Lp{2}\rb{0,\tend;\HM{r_\uu-1}}}}.
\end{align}

For $0\leqslant t\leqslant T$, combining the above estimates, we obtain 
\begin{subequations}
\begin{align}
	\norm{\rb{\uu-\uu_h}\rb{t}}_\LP{\infty}
		\leqslant &Ch \norm{\uu\rb{t}}_\WMP{1}{\infty}
		+C h^{-\nicefrac{d}{2}+r_\uu}\abs{\uu\rb{t}}_\HM{r_\uu}
		\\+&C h^{-\nicefrac{d}{2}+r_\uu-1} e^{C_\uu\nicefrac{\tend}{2}} 
			\max\set{\abs{\uu}_{\Lp{2}\rb{0,\tend;\HM{r_\uu}}},
			\abs{\partial_t\uu}_{\Lp{2}\rb{0,\tend;\HM{r_\uu-1}}}},
\end{align}	
\end{subequations}

from which the claim follows by squaring.

\end{thmProof}

\subsection{Edge-stabilised FEM}	 \label{sec:EdgeStab}

Assuming that the time derivative of the exact solution is only slightly smoother, to be precise, supposing that $\partial_t\uu\in\Lp{2}\rb{0,\tend;\HM{r-\half}}$, the interpolation estimate \eqref{eq:IntResultDtu} can be improved by half an order to 
\begin{equation}
	\norm{\partial_t\uETA\rb{\tau}}_\LP{2}^2  
		\leqslant C_s^2 h^2 \inf_{\vv_h\in\VV_h}\norm{\partial_t\uu\rb{\tau}-\vv_h}_\HM{1}^2
			\leqslant C h^{2r_\uu-1}\abs{\partial_t\uu\rb{\tau}}_\HM{r_\uu-\half}^2.
\end{equation}
Thus, it is estimate \eqref{eq:IntResultConv} that prevents an additional factor $\nu+h$ (instead of $1+h$)  in \eqref{eq:VelDiscConvRate}. The appearance of the unfavourable term $\abs{\uETA\rb{\tau}}_\HM{1}^2$ can be traced back directly to the estimation of the difference of the convective terms in Lemma \ref{lem:ConvTerms} and is therefore due to the nonlinear term in the Navier--Stokes equations. A standard remedy in the case of convection-dominated flows is the introduction of an additional stabilisation term. Let us discuss briefly the potential impact of stabilisation on our method. To this end, with a user-chosen parameter $\gamma\geqslant 0$, we add the term
\begin{equation} \label{eq:EdgeStab}
	s_h\rb{\uu_h,\vv_h} = 
		\sum_{F\in\Fi}\oint_F \gamma h_F^2 \abs{\uu_h\cdot\n_F}^2
		\jmp{\nabla\uu_h}:\jmp{\nabla\vv_h} \ds
\end{equation}
on the left-hand side of the weak formulation on the discrete level \eqref{eq:DiscVarTINS}. Here, $\Fi$ denotes the set of interior edges corresponding to the decomposition $\T$ with outward unit normal vector $\n_F$ and $\jmp{\cdot}$ denotes the standard jump operator. A term such as this is referred to as edge-stabilisation and, because it penalises the jumps of the whole velocity gradient across interior edges, it acts as a stabiliser for problems with dominant convection. Moreover, as shown in \cite{BurmanLinke08}, it is also possible to replace the edge-stabilisation term with a local projection stabilisation term. \\

Now, it is straightforward to show that Theorem \ref{thm:VelDiscErr} still holds true after adding the term $s_h\rb{\uERR,\uERR}$ on the left-hand side, thereby yielding an error of order $k$ for this term as well. An improved error estimate of the difference of convective terms (see Lemma \ref{lem:ConvTerms}) is possible if, for example, the orthogonal $\LTWO$-projector onto $\VV_h$ is used instead of the Stokes projector; cf. \cite{BurmanFernandez07,BurmanLinke08}. Unfortunately, this projector does not act in $\VV_h^\dvg$ and it is not obvious whether it is possible to obtain pressure-robust estimates with this technique. Thus, the existence of an improved pressure- and semi-robust error estimate of order $\mathcal{O}\rb{h^k\sqrt{\nu+h}}$ remains an open question. Nonetheless, we will use convection stabilisation successfully in the form of \eqref{eq:EdgeStab} in the next section concerning numerical examples. 

\section{Numerical experiments}	 \label{sec:NumEx}
In this section, we want to show numerical examples with $\uu\in\Lp{\infty}\rb{0,\tend;\WMP{1}{\infty}}$ for potentially high Reynolds numbers but without consideration of turbulent flows. Recall the space-semidiscrete variational formulation of \eqref{eq:VarTINS} where we replaced the nonlinear term with a generic form $n\rb{\cdot;\cdot,\cdot}$ to be specified later on: 
\begin{subequations} \label{eq:DiscVarTINSrep}
	\begin{empheq}[left=\empheqlbrace]{align} 
	\text{find }\rb{\uu_h,p_h}\colon\rb{0,\tend}\to\VV_h\times\Q_h
	\text{ with }\uu_h\rb{0} =\uu_{0h}	&\text{ s.t. }\forall\rb{\vv_h,q_h}\in\VV_h\times\Q_h\\
		\rb{\partial_t\uu_h,\vv_h} + a\rb{\uu_h,\vv_h} + n\rb{\uu_h;\uu_h,\vv_h} + b\rb{\vv_h,p_h} &=\bra{\ff,\vv_h}  			\\
		  -b\rb{\uu_h,q_h} &= 0 			
	\end{empheq} 
\end{subequations}

In the following, we compare the performance of different standard conforming FEM schemes of order $k$, based on the Taylor--Hood velocity/pressure pair with the Scott--Vogelius pair as a representative of the class of exactly divergence-free and conforming FEM. Inspired by \cite{CharnyiEtAl17}, for the non-divergence-free methods we consider different formulations of the nonlinear term $n\rb{\cdot;\cdot,\cdot}$ which, together with the particular finite element spaces, possess a varying degree of conservation properties and are all equivalent whenever $\nabla\cdot\uu_h=0$; cf. Section \ref{sec:ConservationProperties}. All considered methods are summarised in Table \ref{tab:Methods}. Note that, contrary to all other methods, the EMAC formulation leads to a velocity contribution in the pressure.\\

\begin{table}[h]
\caption{Different finite element spaces and trilinear forms for comparison purposes.}
\label{tab:Methods}
\centering 
\begin{tabular}{cccc} 
\toprule
Name						& $\VV_h/\Q_h$ 		& Trilinear form $n\rb{\uu_h;\uu_h,\vv_h}$	& Abbreviation\\ 
\otoprule
Scott--Vogelius 			& $\sqb{\Pk{k}}^d/\Pdk{k-1}$ 	& $t\rb{\uu_h;\uu_h,\vv_h}$		& $\ger{SV}_k$\\
convective Taylor--Hood 		& $\sqb{\Pk{k}}^d/\Pk{k-1}$ 		& $t\rb{\uu_h;\uu_h,\vv_h}$		& $\ger{convTH}_k$\\ 
skew-symmetric Taylor--Hood & $\sqb{\Pk{k}}^d/\Pk{k-1}$ 		& $t\rb{\uu_h;\uu_h,\vv_h}+\frac{1}{2}\rb{\nabla\cdot\uu_h,\uu_h\cdot\vv_h}$			& $\ger{skewTH}_k$\\  
EMAC Taylor--Hood & $\sqb{\Pk{k}}^d/\Pk{k-1}$ 		& $2\rb{\DD{\uu_h}\uu_h,\vv_h}+\rb{\nabla\cdot\uu_h,\uu_h\cdot\vv_h}$			& $\ger{emacTH}_k$\\
\bottomrule
\end{tabular}
\end{table}

We always guarantee that the starting situation for the SV simulations is in general less favourable compared to the TH simulations by using a coarser mesh for SV such that the total number of DOFs is approximately the same. By doing this, results cannot only be compared in terms of accuracy, but also in terms of efficiency. For numerical simulations, we take advantage of the finite element package COMSOL Multiphysics 5.1. The time stepping is performed using the BDF(2)-scheme and Newton iterations are converged up to a relative residual of $\num{E-8}$. \\

For the purpose of identifying vortical structures in the subsequent flows, we make use of the popular $Q\rb{\uu}$-criterion \cite{HuntEtAl88,JeongHussain95,Haller05}
\begin{equation}
	\set{\x\in\Omega\colon Q\rb{\uu}=\frac{1}{2}\sqb{\abs{\SP{\uu}}^2-\abs{\DD{\uu}}^2}>0},
\end{equation}
that is, we define the neighbourhood of a vortex as the set of points in a flow for which the Euclidean norm of the spin tensor (local rigid body rotation) dominates the deformation tensor (shearing). \\

The first two examples deal with viscous vortex dynamics for $\ff=\zero$ where the approximation of inviscid flows by high Reynolds number viscous flows is considered. Indeed, denoting the inviscid solution of the incompressible Euler equations by $\uu^0$ and the corresponding viscous solution of the incompressible Navier--Stokes equations by $\uu^\nu$, it can be shown \cite[Section 3.1.2]{MajdaBertozzi02} that, for small viscosity,
\begin{align} \label{eq:CompEulerNS}
	\norm{\uu^\nu-\uu^0}_{\Lp{\infty}{\rb{0,\tend;\LP{2}}}}
		\leqslant \nu \norm{\Delta\uu^0}_{\Lp{1}{\rb{0,\tend;\LP{2}}}}
		e^{\abs{\uu^0}_{\Lp{1}\rb{0,\tend;\WMP{1}{\infty}}}}.
\end{align}

\subsection{Gresho-vortex in viscous incompressible flow}	

At first, we consider the Gresho-vortex problem \cite{GreshoChan90,LiskaWendroff03} on $\Omega=\rb{-0.5,0.5}^2$, defined by the evolution of an initial condition that solves the steady inviscid Euler equations. The initial state of the fluid system is originally given in polar coordinates $\rb{r,\phi}$ by means of the velocity 
\begin{align} \label{eq:GreshoInitialVortex}
u_\phi\rb{r,\phi}=
\begin{cases}
5r, \\
2-5r, \\
0,	
\end{cases}
\quad\Rightarrow\quad \uu_0\rb{\x}=\begin{cases}
\rb{-5x_2,5x_1}^\dag,\quad &0\leqslant r\leqslant 0.2,\\
\rb{-\frac{2x_2}{r}+5x_2, \frac{2x_1}{r}-5x_1}^\dag,&0.2\leqslant r\leqslant 0.4,\\
\rb{0,0}^\dag,&0.4\leqslant r,	
\end{cases}
\end{align}
where the radial velocity component $u_r$ vanishes everywhere. The pressure is determined by the zero-mean constraint. This problem is especially interesting because the initial vorticity $\omega=\partial_x u_2-\partial_y u_1$ is discontinuous at $r=\set{0.2,0.4}$ and $\Delta\uu_0\notin\LP{2}{\rb{\Omega}}$, see \eqref{eq:CompEulerNS}. However, with reference to Section \ref{sec:ErrorEstimatesVel}, note that $\uu_0\in\WMP{1}{\infty}{\rb{\Omega}}$ with $\abs{\uu_0}_\WMP{1}{\infty}=5$.\\

Now, our numerical experiment consists in taking this vortex \eqref{eq:GreshoInitialVortex} as the initial condition of the simulation, that is $\uu_h\rb{0}=\uu_0$, and letting it evolve within a viscous incompressible flow. Then, depending on the viscosity $\nu>0$, more or less viscous dissipation causes an alteration of the vortex flow and may give rise to numerical instabilities. This is a different scenario than the one considered in \cite{CharnyiEtAl17,LiskaWendroff03} since we are deliberately including viscous effects and also observe evolutions over longer time periods.  \\
 
As regards the mesh construction, we always begin with a decomposition into $N\times N$ quadratic elements which is transformed into a simplicial mesh by inserting the diagonals and, for the SV-FEM, is then barycentre-refined. The meshes for the TH-FEM are referred to as $N/$ whereas the barycentre-refined meshes are denoted by $N/b$. We restrict the BDF(2) maximum time step to $\Delta t_\maxrm=0.005$. \\

\begin{figure}[h]
\centering
\includegraphics[width=0.975\textwidth]{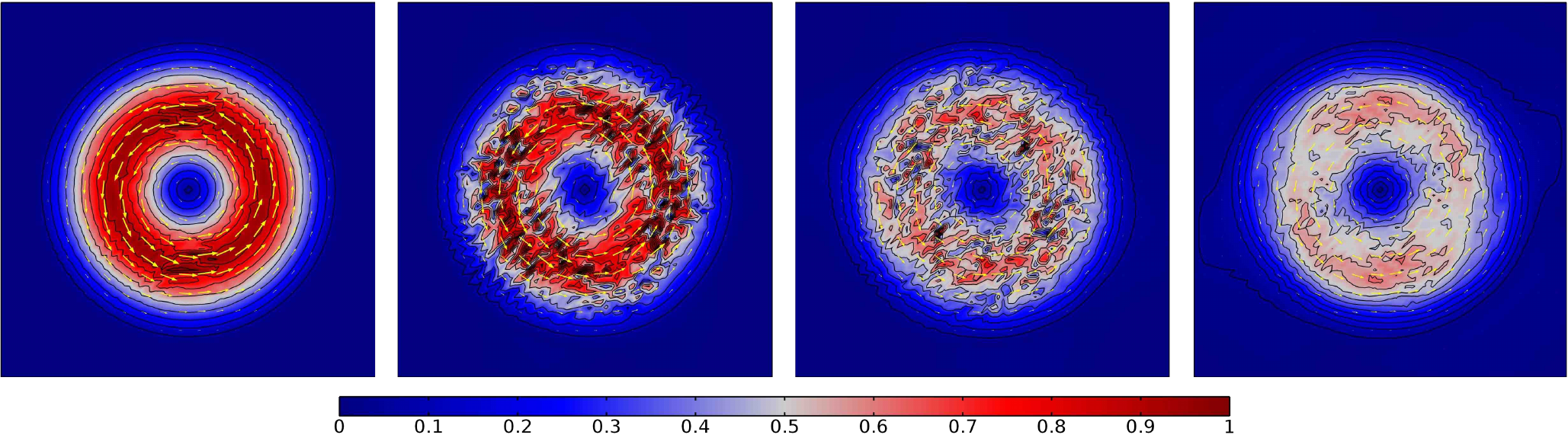}
\caption{$\abs{\uu_h}$-evolution of Gresho-vortex for $\nu=\num{7E-5}$ with $\ger{skewTH}_2$ on $70/$ at $t\in\set{1,2,3,4}$.}
\label{fig:skewTH2}
\end{figure}

In Figure \ref{fig:skewTH2}, we see the evolution of the velocity magnitude $\abs{\uu_h}$ for a $\ger{skewTH}_2$ method with $\nu=\num{7E-5}$ on a $70/$ mesh (\num{44803} DOFs). It can be observed that, as time proceeds, the vortex looses its shape and is dissipated quickly and in a rather fragmented manner. As we will see later on, this decay happens much too fast and is due to the unsuitability of the skew-symmetric formulation for such a problem. Note that, even though this formulation conserves kinetic energy, this property is obviously not strong enough in the context of vortex dynamics. Moreover, conducted tests with higher-order finite element pairs revealed no improvement for the skew-symmetric formulation, which is possibly due to the lack of smoothness of $\uu_0$. \\

\begin{figure}[h]
\centering
\subfigure[$\ger{convTH}_2$, $70/$, $t=4$]{\includegraphics[width=0.23\textwidth]{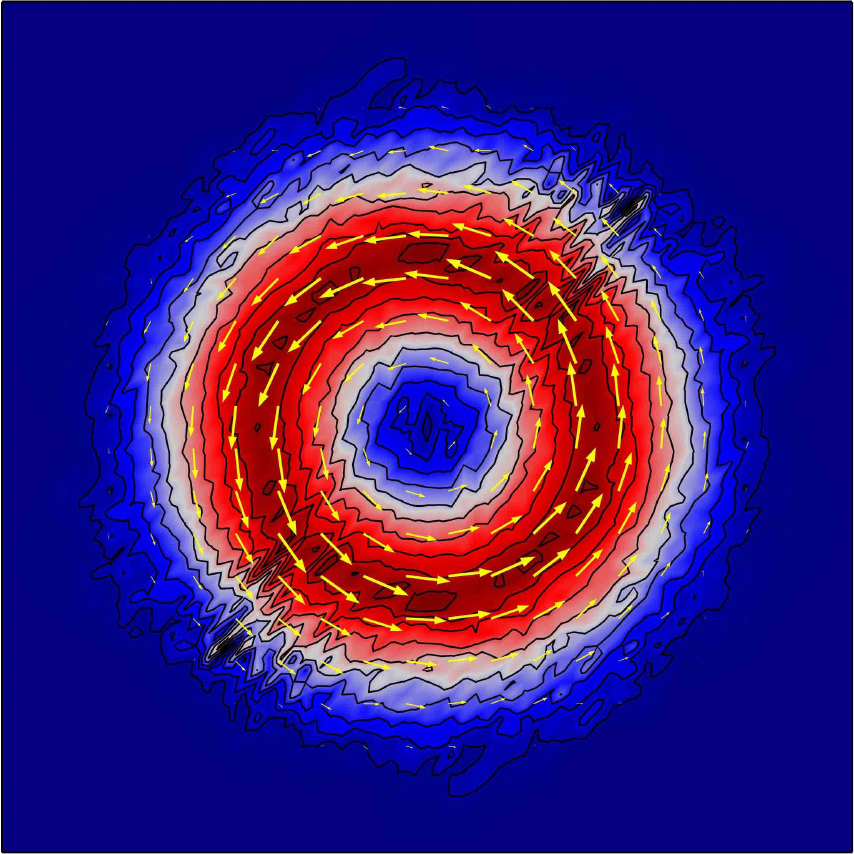}} \goodgap
\subfigure[$\ger{emacTH}_2$, $70/$, $t=4$]{\includegraphics[width=0.23\textwidth]{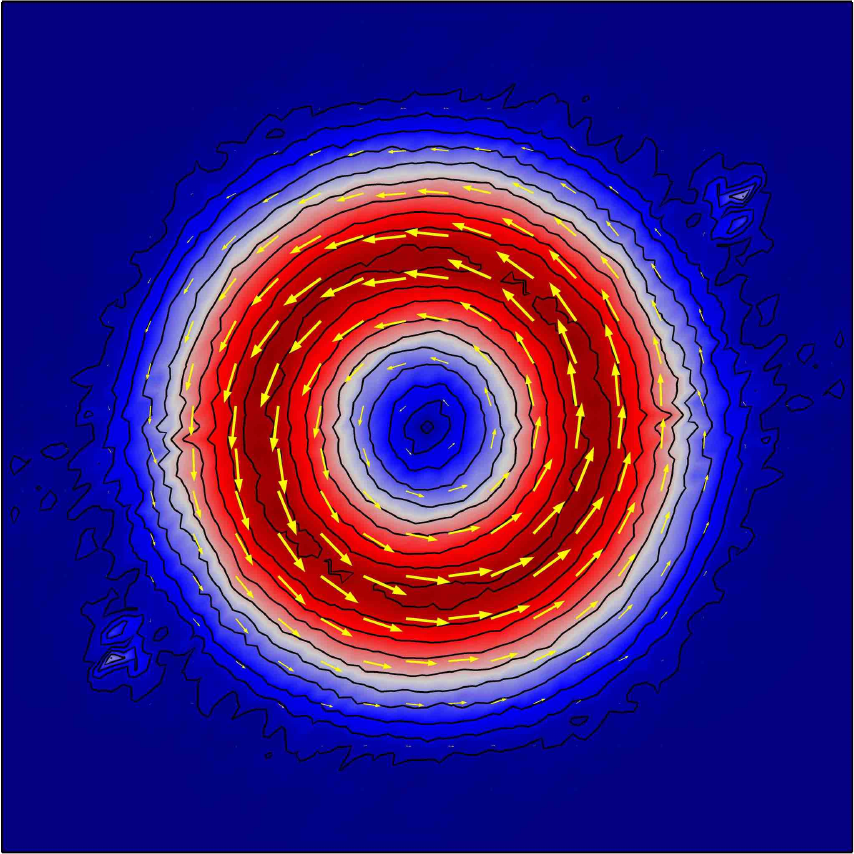}} \goodgap
\subfigure[$\ger{convTH}_3$, $45/$, $t=2.7$]{\includegraphics[width=0.23\textwidth]{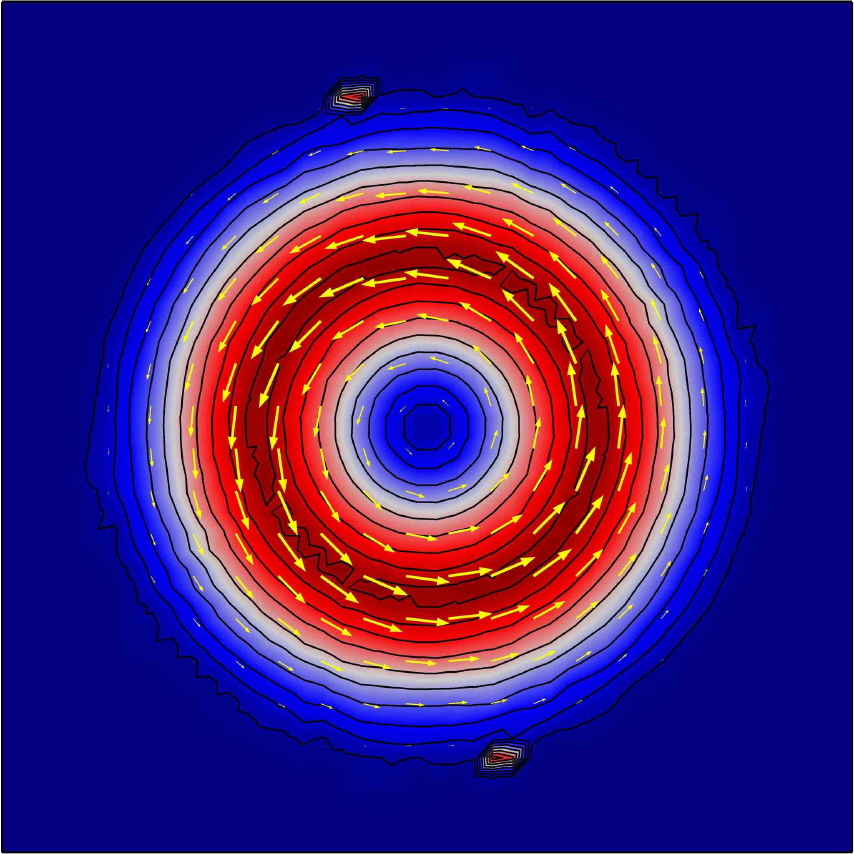}} \goodgap
\subfigure[$\ger{emacTH}_3$, $45/$, $t=0.8$]{\includegraphics[width=0.23\textwidth]{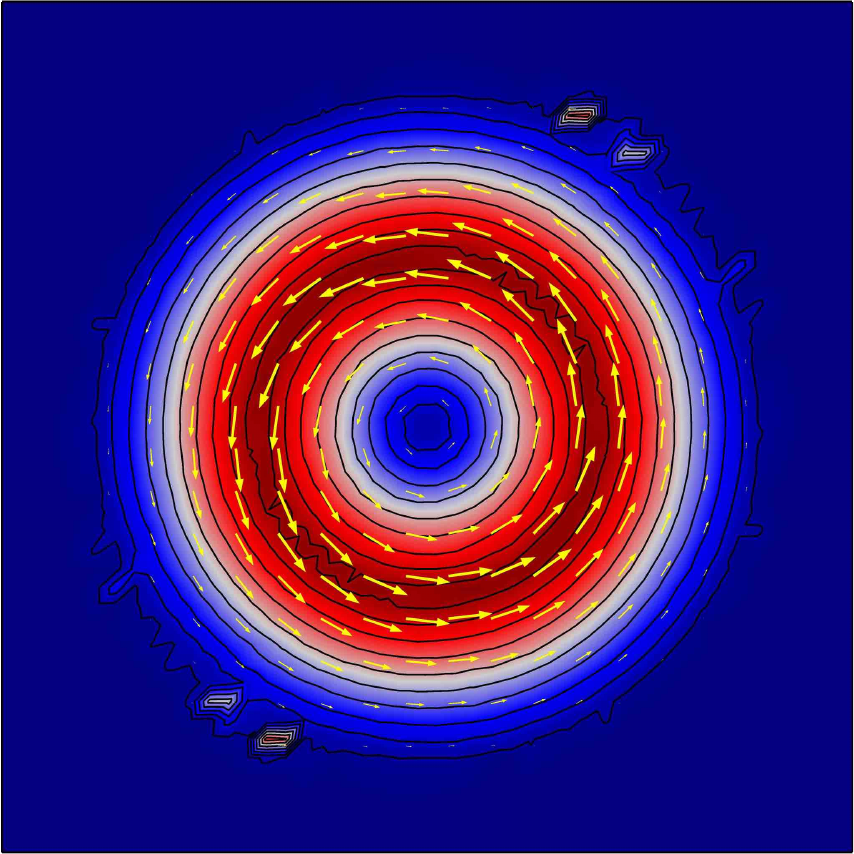}} 
\caption{Comparison of $\abs{\uu_h}$ of Gresho-vortex for $\nu=\num{4E-6}$. Colour bar identical to other figures.}
\label{fig:ConvEMAC}
\end{figure}

However, if we use the convective and EMAC formulation of the nonlinear term for the same problem, the vortex is still rather distinct at $t=4$ and therefore, in Figure \ref{fig:ConvEMAC}, in order to be able to observe unusual behaviour, we regard the problem for a smaller viscosity of $\nu=\num{4E-6}$. In the first two images, we compare second-order $\Pk{2}/\Pk{1}$ TH-FEM with different nonlinear terms at the same time $t=4$. Obviously, both methods are better suited to this problem than the skew-symmetric formulation. However, the vortex from $\ger{convTH}_2$ is also not truly stable and, for $\ger{emacTH}_2$, we see the outlines of artefacts arising outside of the vortex, which oftentimes lead to divergence of the nonlinear solver as time proceeds.  \\

As can be seen in the third and fourth image in Figure \ref{fig:ConvEMAC}, the use of a higher-order $\Pk{3}/\Pk{2}$ pair does not resolve the problem either. In order to produce comparable situations we use a coarser $45/$ mesh (\num{45273} DOFs) for the third-order methods. While the vortex itself is resolved well in both cases and especially its boundary seems to be rather sharp, both convective and EMAC formulations yield the already mentioned artefacts. Unfortunately, the nonlinear solvers failed before reaching $t=4$ at the depicted time instances, in each case. An improvement of both methods can of course be obtained by refining the mesh which yields more DOFs, thereby making these standard methods computationally inefficient. \\

\begin{figure}[h]
\centering
\includegraphics[width=0.975\textwidth]{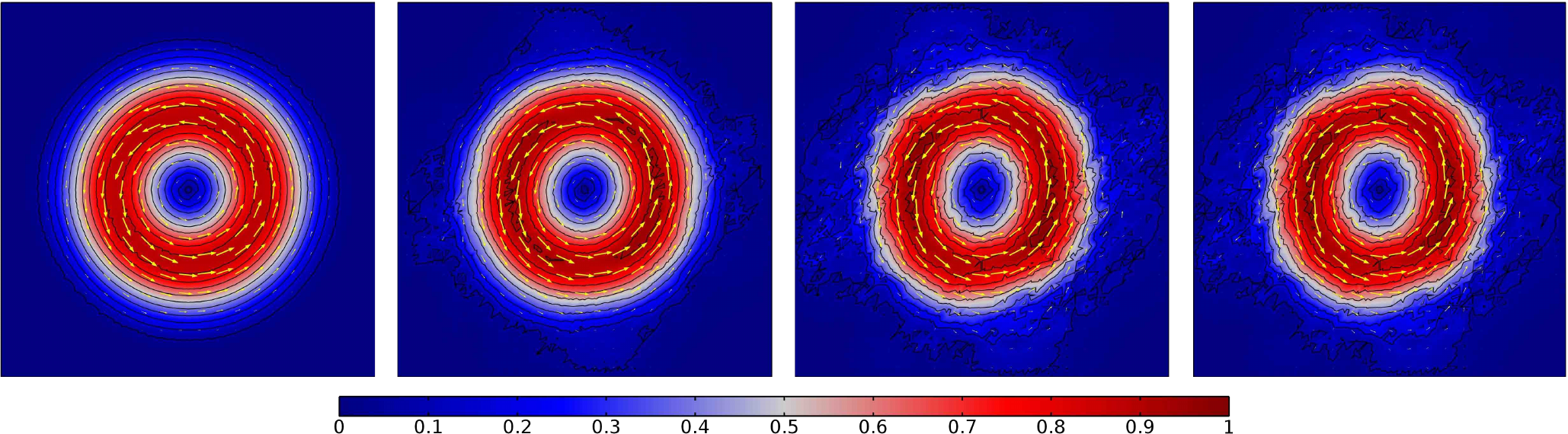}
\caption{$\abs{\uu_h}$ of Gresho-vortex at $t=4$ for $\nu\in\set{\num{7E-5},\num{4E-6},\num{E-8},\num{0}}$ with $\ger{SV}_2$ on $32/b$.}
\label{fig:SV2}
\end{figure}

Finally turning our attention to the exactly divergence-free SV-FEM, in Figure \ref{fig:SV2} we compare the velocity magnitude at $t=4$ for the $\ger{SV}_2$ method on a $32/b$ mesh (\num{43266} DOFs) for different viscosities. The first image represents the SV solution corresponding to the viscosity considered for the $\ger{skewTH}_2$ case and the second image can be compared to the $\ger{convTH}_2$ and $\ger{emacTH}_2$ results. It can be observed that the exactly divergence-free method is clearly superior in both cases, as neither considerable smearing or dissipation of the vortex occurs nor any artefacts appear to incapacitate the method. In addition, we can even decrease the viscosity for $\ger{SV}_2$ until reaching $\nu=0$ without losing the vortical characteristic. Interestingly, it is not possible to see any significant differences between $\nu=0$ and $\nu=\num{E-8}$. Of course, diffusive effects result in non-zero velocity outside the vortex, which is also slightly distorted; but nevertheless, we believe that these are very convincing and satisfactory results for the SV-FEM. \\

\begin{figure}[h]
\centering
\includegraphics[width=0.975\textwidth]{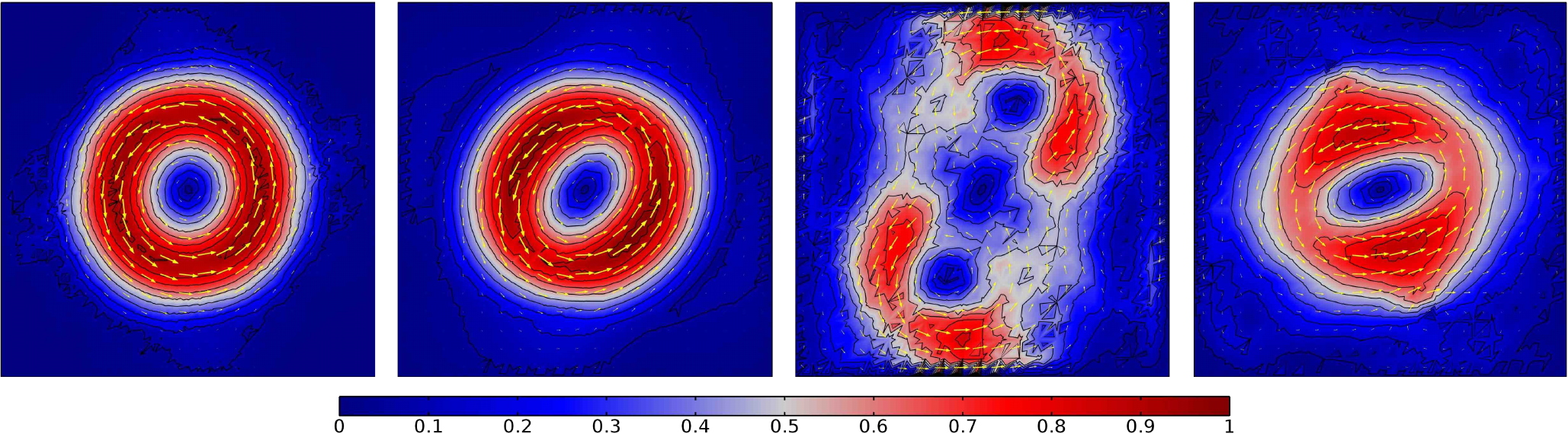}
\caption{$\abs{\uu_h}$-evolution of Gresho-vortex for $\nu=\num{4E-6}$ with $\ger{SV}_4$ on $20/b$ at $t\in\set{4,8,16,32}$.}
\label{fig:SV3}
\end{figure}

Lastly, we regard the long-term behaviour of the SV-FEM. To this end, the evolution of $\abs{\uu_h}$ for a fourth-order $\ger{SV}_4$ with viscosity $\nu=\num{4E-6}$ on a $20/b$ mesh (\num{62723} DOFs) is shown in Figure \ref{fig:SV3}. First of all, we note that, in contrast to Taylor--Hood, higher-order finite element pairs work fine for Scott--Vogelius. Secondly, we see that the vortex is stable up to $t=8$, which is remarkable. Going beyond this time instance, the vortex dissolves as can be seen at $t=16$. However, energy and momentum of the vortex are not dissipated as, for example, is clearly the case for $\ger{skewTH}_2$. Indeed, a vortical structure somehow reassembles until $t=32$ and it is surprising that the SV-FEM is able to reliably produce results for such a long period of time. 

\subsection{Dynamics of planar lattice flow}	
\begin{figure}[h]
\centering
\includegraphics[width=0.975\textwidth]{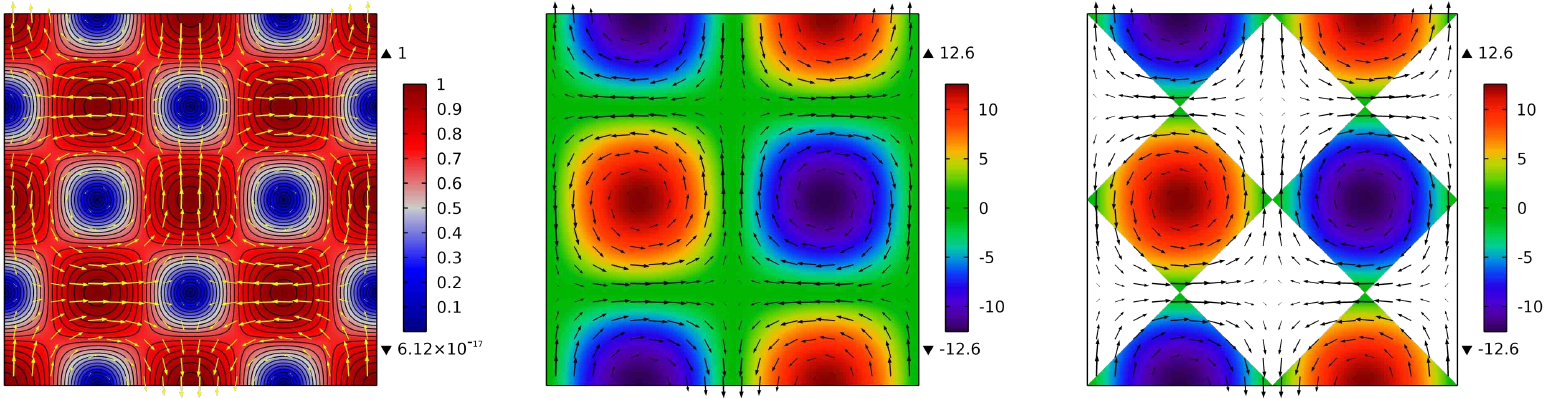}
\caption{Initial state of planar lattice flow: Velocity magnitude $\abs{\uu_0}$, vorticity $\omega_0$ and $Q\rb{\uu_0}$-criterion coloured with $\omega_0$. The arrows denote the velocity field and the triangles denote the maximum and minimum value attained.}
\label{fig:ICLatticeFlow}
\end{figure}

In this section, we again consider the evolution of an initial velocity, which solves the stationary incompressible Euler equation, in a viscous incompressible Navier--Stokes flow. However, in this example, the following exact solution $\uu$ for $\nu\geqslant 0$ and corresponding initial condition $\uu_0$ is known:
\begin{equation} \label{eq:PlanarExactSol}
\uu_0\rb{\x}=\begin{pmatrix}
\sin\rb{2\pi x_1}\sin\rb{2\pi x_2}\\
\cos\rb{2\pi x_1}\cos\rb{2\pi x_2}	
\end{pmatrix},\quad
\uu\rb{t,\x} = \uu_0\rb{\x}e^{-8\pi^2 \nu t },\quad
\x\in\Omega=\rb{0,1}^2
\end{equation}
The initial velocity field induces a flow structure called `planar lattice flow' \cite{Bertozzi88} which, due to its saddle point character, is `dynamically unstable so that small perturbations result in a very chaotic motion' \cite{MajdaBertozzi02}. Here, we impose periodic boundary conditions on the vertical and horizontal walls of $\partial\Omega$, respectively, and the zero-mean condition on $p_h$. Extending the field 1-periodically to $\R^2$, we obtain $\uu_0\in\CK{\infty}{\rb{\R^2}}$ and $\abs{\uu}_{\Lp{\infty}\rb{0,\tend;\WMP{1}{\infty}}}=\abs{\uu_0}_\WMP{1}{\infty}=2\pi$ for any $\tend>0$. In Figure \ref{fig:ICLatticeFlow}, the velocity magnitude $\abs{\uu_0}$, the initial vorticity $\omega_0$ and the $Q\rb{\uu_0}$-criterion coloured with the vorticity are shown. Similarly as for the Gresho-vortex, the dynamical behaviour of such a system is observed when it is exposed to a viscous flow. \\

Contrary to the previous subsection, we now consider unstructured Delaunay meshes where the separatrices of the flow (here, zero level set of $\omega_0$) are deliberately not aligned with the mesh. The BDF(2) is restricted to a slightly larger maximum time step of $\Delta t_\maxrm=0.01$. For the sake of brevity, only a fixed viscosity $\nu=\num{4E-6}$ is considered and the Scott--Vogelius FEM is compared directly to the EMAC Taylor--Hood FEM. We briefly demonstrate that the other TH formulations are practically useless for this problem.\\

\begin{figure}[h]
\centering
\includegraphics[width=1\textwidth]{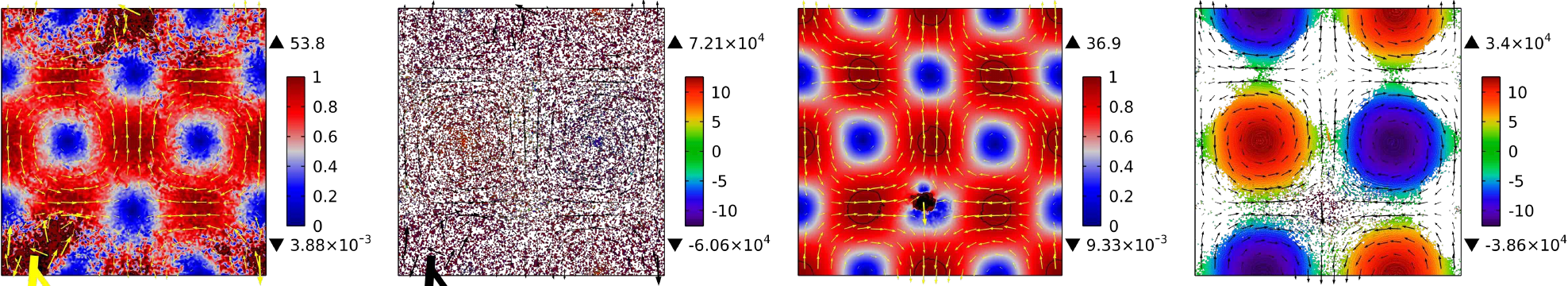}
\caption{Poor numerical solutions for planar lattice flow: Velocity magnitude and $Q\rb{\uu_h}$-criterion coloured with $\omega$ for $\ger{skewTH}_2$ at $t=0.47$ (first and second image) and $\ger{convTH}_2$ at $t=6.4$ (third and forth image) using \num{194395} DOFs.}
\label{fig:SkewConvLatticeFlow}
\end{figure}

\begin{figure}[h]
\centering
\subfigure{\includegraphics[width=0.49\textwidth]{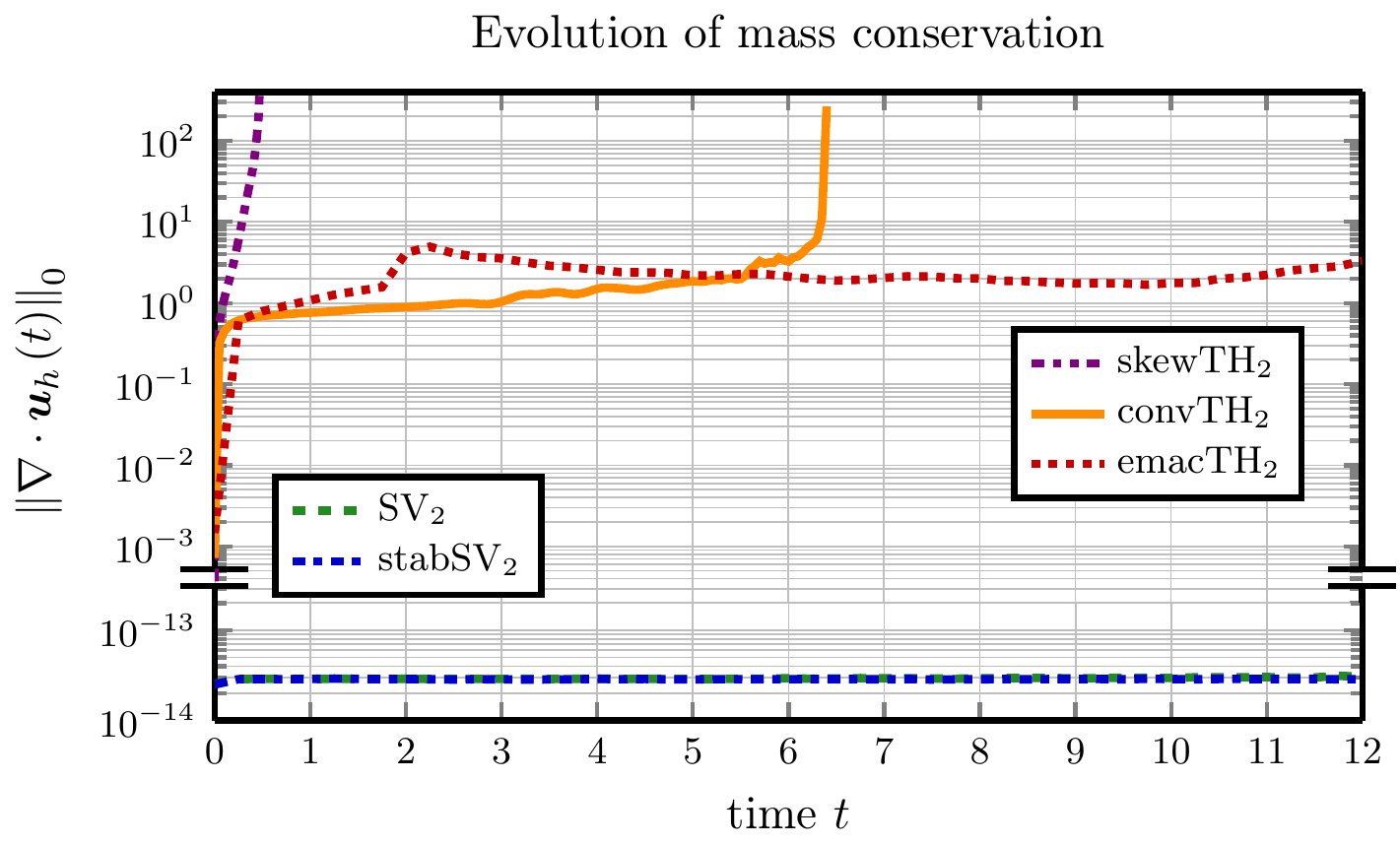}} \goodgap
\subfigure{\includegraphics[width=0.49\textwidth]{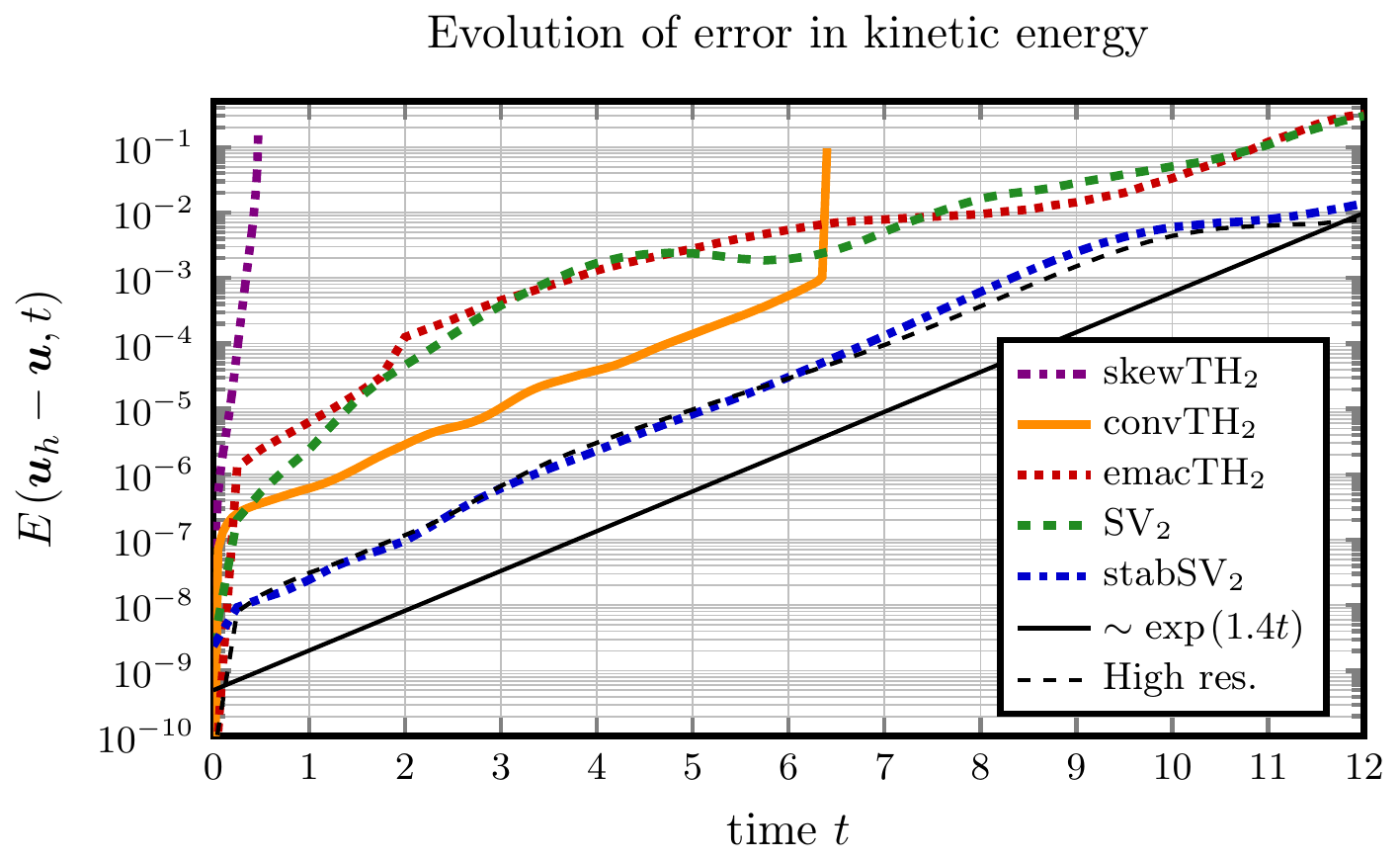}} \\
\subfigure{\includegraphics[width=0.49\textwidth]{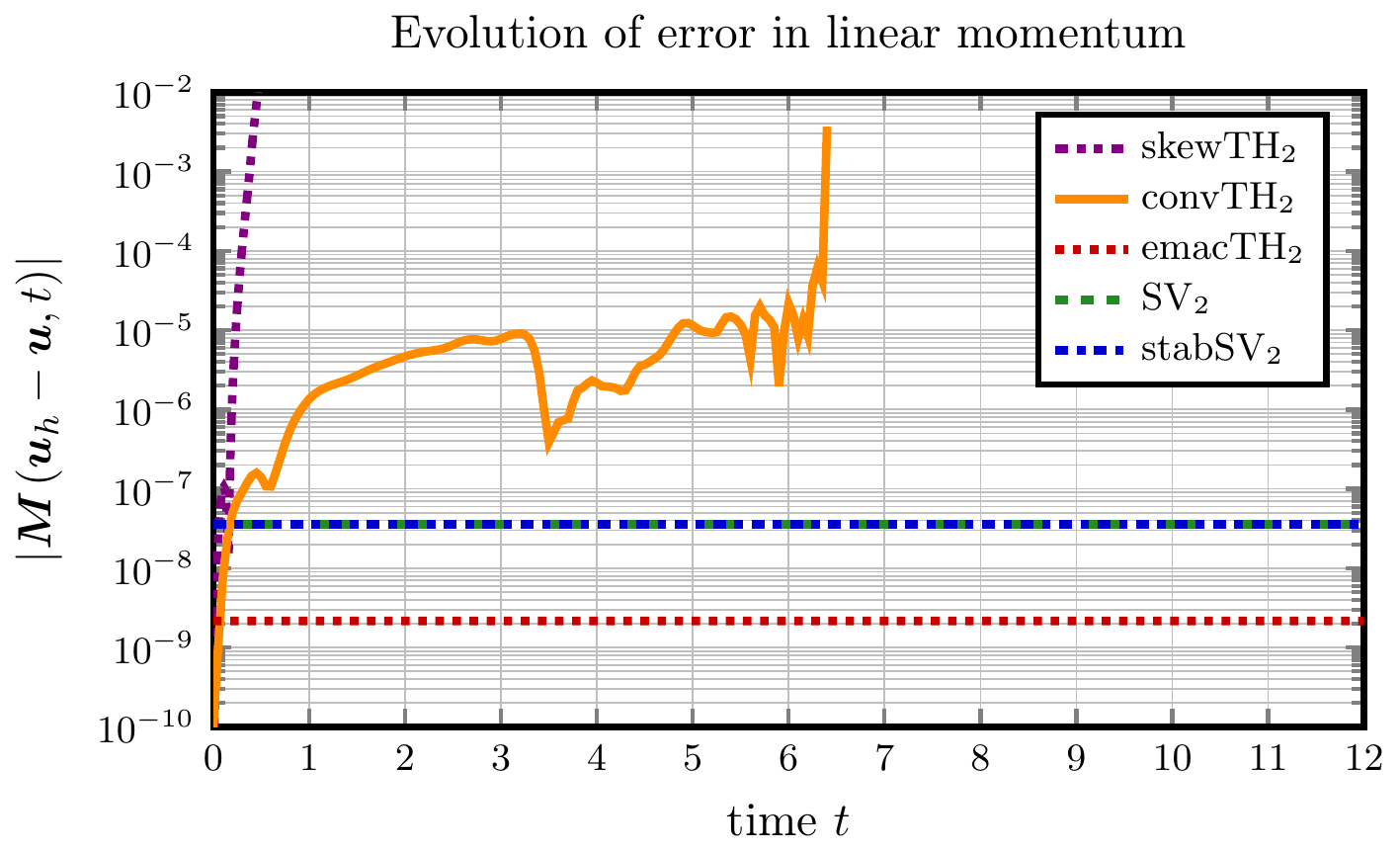}} \goodgap
\subfigure{\includegraphics[width=0.49\textwidth]{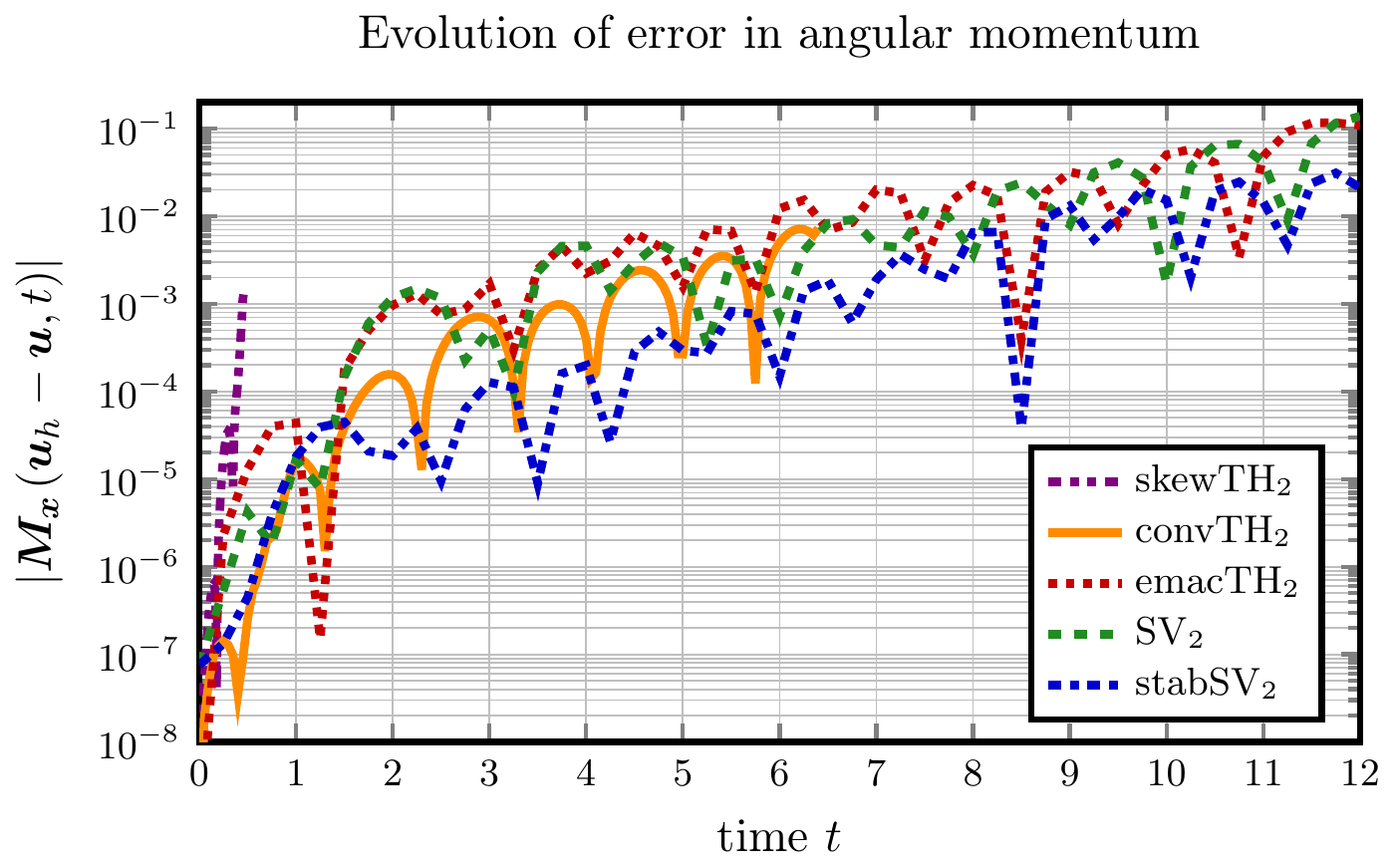}} 
\caption{Comparison of the errors w.r.t.\ the conservation properties (Section \ref{sec:ConservationProperties}) for all different methods with $\nu=\num{4E-6}$.}
\label{fig:Conservation}
\end{figure}

In order to show the unsuitability of skew-symmetric and convective formulations of the nonlinear term, we take a fine mesh with $h=\num{1.17E-2}$ (\num{194395} DOFs) and perform a simulation for the planar lattice flow. Both methods fail to converge and the results, shown at the particular last possible time instance, can be seen in Figure \ref{fig:SkewConvLatticeFlow}. Whilst the skew-symmetric formulation fails at an early stage and yields an obviously poor solution, the convective formulation remains stable longer and the solution in principal seems to be better. However, comparing the maximum and minimum values in the $Q\rb{\uu_h}$-criterion, we observe that the vorticity blows up locally in comparison to the initial vorticity $\omega_0$. We shall see that this phenomenon is not as pronounced for either the EMAC or SV method. Note that, in Figure \ref{fig:Conservation}, $\ger{skewTH}_2$ and $\ger{convTH}_2$ are also taken into account in a quantitative comparison.\\

\begin{figure}[h]
\centering
\includegraphics[width=0.975\textwidth]{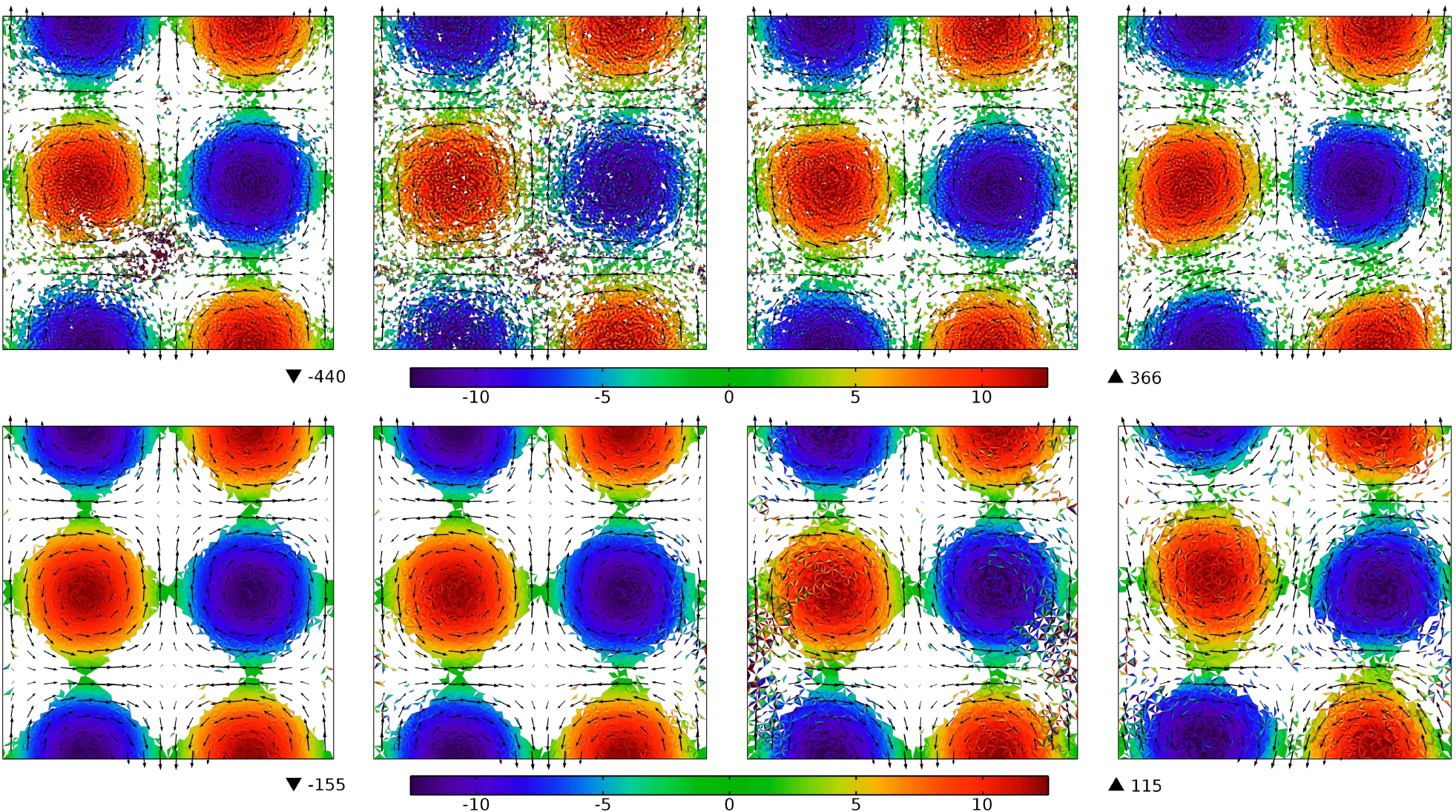}
\caption{Evolution of $Q\rb{\uu_h}$-criterion coloured with $\omega$ at $t\in\set{2,3,6,9}$ (columns). Comparison: $\ger{emacTH}_2$ using \num{61808} DOFs (upper row) and $\ger{SV}_2$ using \num{59766} DOFs (lower row).}
\label{fig:Q-comparison}
\end{figure}

In Figure \ref{fig:Q-comparison}, we compare the performance of the EMAC-TH and SV method by means of the evolution of the $Q\rb{\uu_h}$-criterion and the maximum values of the corresponding vorticity. The $\ger{emacTH}_2$ solution is computed on a mesh with $h=\num{2.14E-2}$ (\num{61808} DOFs) and $\ger{SV}_2$ uses $h=\num{4.47E-2}$ (\num{59766} DOFs). Firstly, we note that, despite using much coarser meshes, EMAC and SV do not fail to converge. However, regarding $\ger{emacTH}_2$ in the upper row, we observe that the initially diamond-shaped form of the $Q\rb{\uu_h}$-criterion breaks apart at $t=2$ and at $t=9$, we see a circular structure. \\

Compared with this, the $Q\rb{\uu_h}$-criterion of $\ger{SV}_2$ remains more diamond-shaped and even at $t=9$, the structure is not as circular as for the EMAC method. Regarding the values of the vorticity we observe that, although the SV method yields an increased vorticity compared to the initial vorticity, its minimum and maximum at least have the same absolute value. The vorticity of the EMAC method, on the other hand, is significantly increased and the absolute values do not coincide. Based on these insights we conclude that the SV-FEM is superior to all TH-type FEM for this kind of planar lattice flow. However, for TH-FEM, the EMAC formulation clearly outperforms the convective and skew-symmetric formulations. \\

\begin{figure}[h]
\centering
\includegraphics[width=0.975\textwidth]{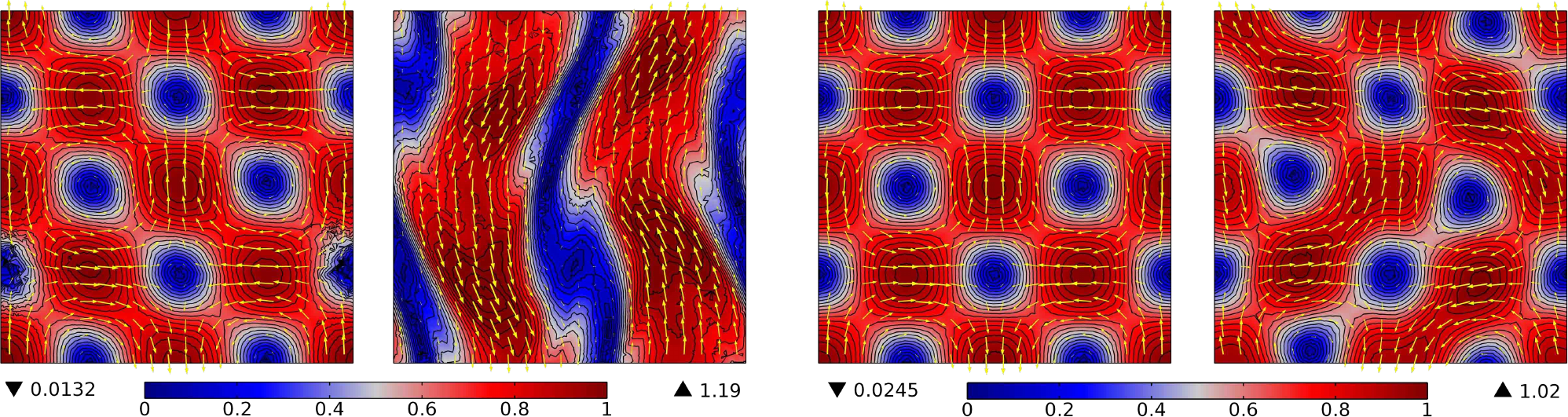}
\caption{Velocity magnitude at $t\in\set{6,12}$ with $\ger{SV}_2$ (\num{59766} DOFs). Comparison: Unstabilised method (left group) and edge-stabilised method with $\gamma=\num{E-1}$ (right group).}
\label{fig:Stab-comparison}
\end{figure}

With regard to Section \ref{sec:EdgeStab}, let us briefly demonstrate the impact of convection stabilisation. Therefore, edge-stabilisation is now added to our method, referred to as $\ger{stabSV}_2$. In Figure \ref{fig:Stab-comparison}, we compare the performance of the unstabilised and edge-stabilised ($\gamma=\num{E-1}$) SV method by means of the evolution of $\abs{\uu_h}$. On the left-hand side, we observe that, as time proceeds, the initial lattice structure of the flow begins to `vibrate' until, eventually, it breaks apart and a chaotic motion develops for the unstabilised $\ger{SV}_2$ method. The flow at $t=12$ does not have any resemblance to a vortical structure anymore. Furthermore, at this time, we get `overshoots' of about \SI{19}{\percent} as the velocity magnitude is overpredicted by the unstabilised method.\\

This instability can be treated with a stabilising term of the form \eqref{eq:EdgeStab}. Indeed, the right-hand side of Figure \ref{fig:Stab-comparison} shows the impact of the addition of such a term with $\gamma=\num{E-1}$ to the $\ger{SV}_2$ method. It can be seen that, although the same kind of `vibration' occurs, the edge-stabilisation is able to counteract the impending collapse of the vortex structure. Thus, at $t=12$ no chaotic behaviour is visible and the overshoots are reduced to \SI{2}{\percent}. However, the vortices at $t=12$ are about to break loose from the lattice, which ultimately leads to them merging; therefore, edge-stabilisation can only postpone, not prevent, the breakdown of the initial planar lattice flow. \\

\begin{figure}[h]
\centering
\includegraphics[width=0.975\textwidth]{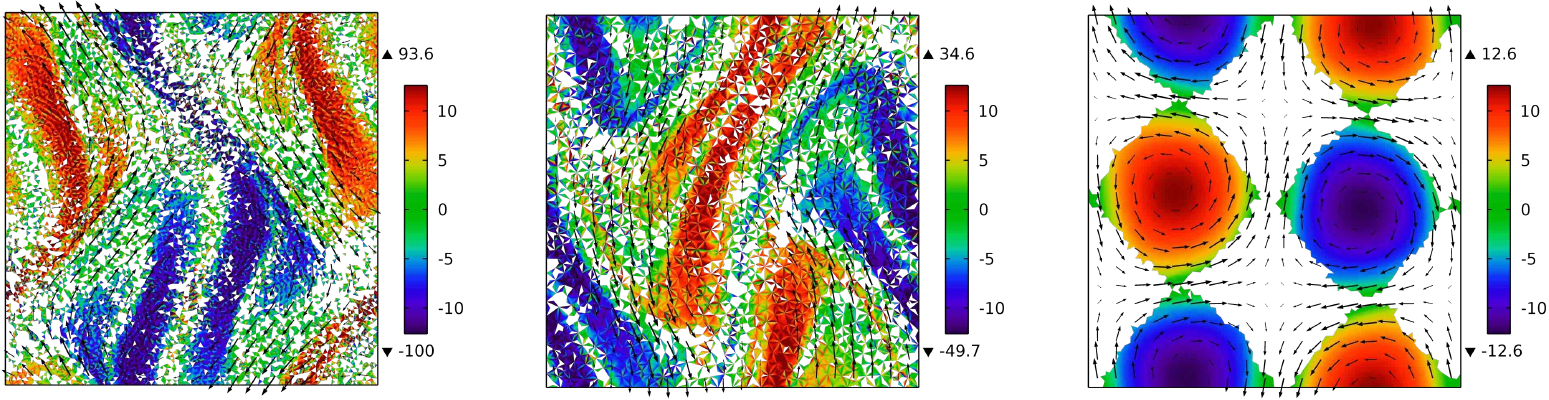}
\caption{$Q\rb{\uu_h}$-criterion coloured with $\omega$ at $t=12$. Comparison (left to right): $\ger{emacTH}_2$ (\num{61808} DOFs), unstabilised $\ger{SV}_2$ and edge-stabilised $\ger{SV}_2$ with $\gamma=\num{E-1}$ (both \num{59766} DOFs).}
\label{fig:Q-Stab-comparison}
\end{figure}

Finally, we compare the $Q\rb{\uu_h}$-criterion and vorticity at $t=12$ for $\ger{emacTH}_2$, unstabilised $\ger{SV}_2$ and $\ger{stabSV}_2$ in Figure \ref{fig:Q-Stab-comparison} and observe that only the stabilised method p†reserves the initial diamond shape. Furthermore, the maximum and minimum values $\pm 12.6$ of the vorticity still exactly coincide with the initial values, see Figure \ref{fig:ICLatticeFlow}, and the $Q\rb{\uu_h}$-criterion does not indicate an accumulation of diffusive vortical structures in between the diamonds. In contrast, the vorticity of both unstabilised $\ger{SV}_2$ and $\ger{emacTH}_2$ increases and the $Q\rb{\uu_h}$-criterion shows a chaotic structure. Therefore, we have proven that stabilisation of the Scott--Vogelius method clearly and visibly improves the quality of the approximate solution in practical applications. \\

The last observations are made with regard to Section \ref{sec:ConservationProperties}. Since we know the exact solution \eqref{eq:PlanarExactSol} to the problem, we can also compute the evolution of its kinetic energy, linear momentum and angular momentum (of course, $\nabla\cdot\uu=0$). Thus, in Figure \ref{fig:Conservation}, it is observable how the errors of the different FEMs w.r.t.\ those conservation properties evolve over time. Note that the exact solution has vanishing linear and angular momentum, i.e.\ $\MM\rb{\uu,t}=\MM_\x\rb{\uu,t}=\zero$ for all $t\in\sqb{0,\tend}$, which does not hold for the FEM approximation. Corresponding to the already shown failure of $\ger{skewTH}_2$, it can be seen that divergence, kinetic energy and linear momentum of $\uu_h$ explode after a short time ($t=0.47$). For the $\ger{convTH}_2$ method the same blow-up is apparent, but occurs after a considerable longer time ($t=6.4$). At this point we note that, as long as $\ger{convTH}_2$ gives results, the error in the kinetic energy is relatively low. However, the simulations with the skew-symmetric and convective nonlinear term are based on a very fine mesh which yields about three times more DOFs than the other formulations. As expected, the conservation properties of $\ger{emacTH}_2$, $\ger{SV}_2$ and $\ger{stabSV}_2$ are clearly superior and do not demonstrate blow-up in either divergence, kinetic energy or linear momentum. The divergence error of the EMAC method is bounded in time and, by construction, the Scott--Vogelius methods have a divergence close to machine precision. Concerning the linear momentum, $\ger{SV}_2$ and $\ger{stabSV}_2$ yield the same error and $\ger{emacTH}_2$ seems to be better; but the overall magnitude is so small anyway that this does not seem to make a difference in the results. The error plot of the angular momentum shows comparable results for the unstabilised SV and the EMAC-TH whereas the edge-stabilised SV-FEM, in general, is slightly less error-prone.\\

Perhaps the most interesting comparison can be observed for the error in the kinetic energy. Here, we added the black solid line which represents a function proportional to $\exp{\rb{1.4 t}}$ with the purpose of verifying experimentally that the error in the kinetic energy increases exponentially in time. The corresponding theoretical prediction manifests in the Gronwall exponent of Theorem \ref{thm:VelDiscErr}. Obviously, the errors of $\ger{emacTH}_2$, $\ger{SV}_2$ and $\ger{stabSV}_2$ show this exponential behaviour also in practice which is a very interesting observation. In addition, it becomes quantitatively apparent that the edge-stabilised method is far more accurate than the pure Galerkin method. Lastly, the dashed black line results from a stabilised SV simulation with high resolution in both space and time (\num{226656} DOFs and $\Delta t_\maxrm=\num{E-3}$). Since the slope is similar to $\ger{stabSV}_2$ we infer that the exponential growth of the error in the kinetic energy is invariant against mesh refinement and temporal discretisation. Moreover, the total error of the high resolution simulation did not decrease in comparison to the $\ger{stabSV}_2$ method. Therefore, we conclude that the exponential term in Theorem \ref{thm:VelDiscErr} dominates the total error and its dynamics. Simultaneously, this problem serves as a warning example, since neither mesh refinement nor decreasing the time step size significantly improves the quality of the approximation.

\subsection{Laminar Blasius boundary layer}	
Let us consider the laminar Blasius boundary layer, that is, the viscous two-dimensional and stationary flow which forms around a semi-infinite flat plate parallel to the incident free-stream. Denoting the velocity of the free-stream by $u_\infty$, the attached laminar boundary layer flow is known to be a self-similar solution. The dimensionless similarity variable is given by $\eta=x_2\sqrt{u_\infty/\rb{\nu x_1}}$ and the velocity, in the vicinity of the plate, can be written as $\uu=u_\infty f^\prime\rb{\eta}$. Here, $f$ solves a third-order ordinary differential equation which can be derived from Prandtl's boundary layer equations \cite{Tritton88,SchlichtingGersten00,Durst08} and we note that $\abs{\uu}_\WMP{1}{\infty}\sim\nu^{-\half}$.\\

\begin{figure}[h]
\centering
\includegraphics[width=0.6\textwidth]{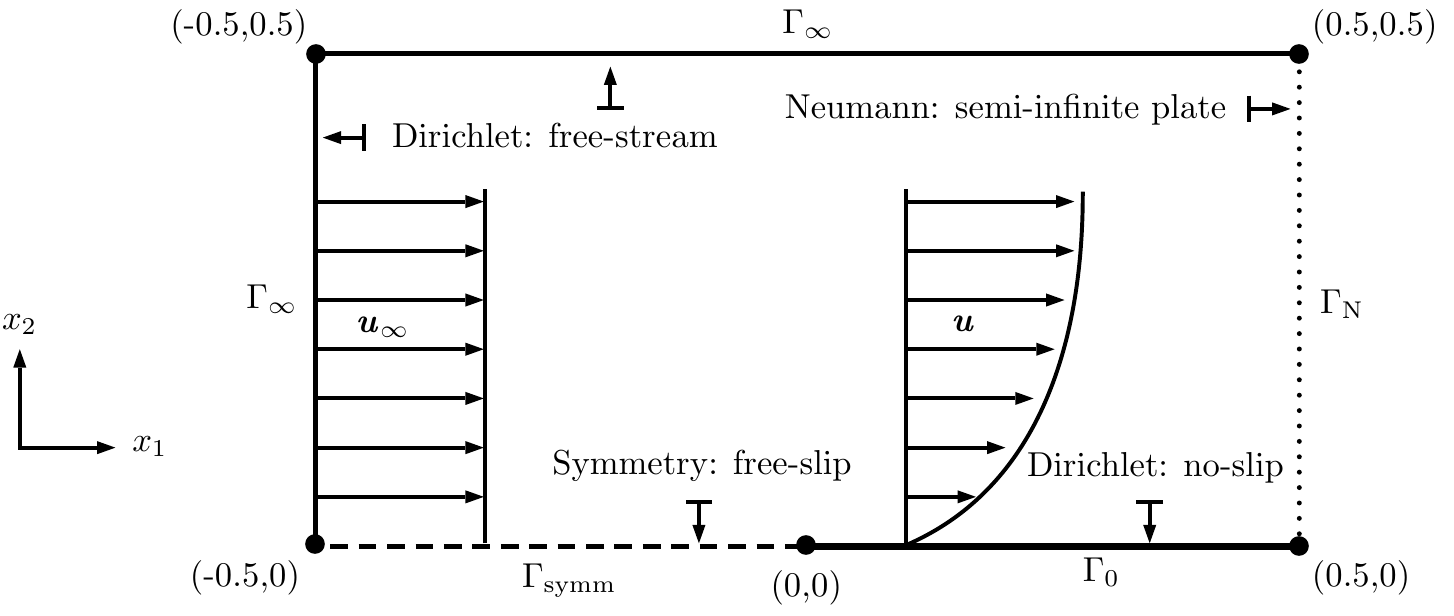}
\caption{Schematic representation of the setup for the laminar Blasius boundary layer problem.}
\label{fig:BlasiusSchematic}
\end{figure}

We now want to see whether the Scott--Vogelius method gives reliable numerical results for such a problem. Comparable computations for Taylor--Hood elements can be found in \cite{ArndtEtAl15,DallmannEtAl16}. In this work, we restrict ourselves to the domain $\Omega=\rb{-0.5,0.5}\times\rb{0,0.5}$ which is actually the upper half of the domain considered in the aforementioned publications. A schematic representation of the setup can be taken from Figure \ref{fig:BlasiusSchematic} and, as can be seen there, the whole problem is determined by its boundary conditions on $\partial\Omega=\Gamma$. Using the free-stream, we impose the Dirichlet condition $\uu=\rb{u_\infty,0}^\dag=\rb{1,0}^\dag$ on $\Gamma_\infty$ and on $\Gamma_0$ we prescribe the no-slip condition $\uu=\zero$ which represents the flat plate. On $\Gamma_\mathrm{N}$ we artificially truncate the plate and use a homogeneous Neumann (do-nothing) condition to mimic a semi-infinite plate. For the symmetry plane $\Gamma_\mathrm{symm}$ we follow \cite{GreshoSani00,Gunzburger89} and use the decomposition $\uu=\rb{\uu\cdot\n}\n+\n\times\rb{\uu\times\n}$ to define the projection of $\uu$ onto the tangent plane by $\uu_\TAU=\uu-\rb{\uu\cdot\n}\n$. Hence, we impose the free-slip condition
\begin{align}
	\uu\cdot\n = 0,\quad  
	2\nu\DD{\uu}\n-\sqb{2\nu\DD{\uu}\n\cdot\n}\n =
	\nabla\uu_\TAU\cdot\n = \frac{\partial \uu_\TAU}{\partial\n}= \zero,
\end{align}	
which ensures vanishing normal velocity (no flow across the boundary) and vanishing shear stress (no viscous stress in tangential direction) on $\Gamma_\mathrm{symm}$. Note that the simplification in the shear stress condition is only valid for planar boundaries \cite{GreshoSani00}. The pressure is already determined uniquely by the do-nothing condition. \\

As already mentioned, our intent is to consider laminar boundary layers and, therefore, our simulations must remain below a certain critical Reynolds number. The interesting Reynolds number is formed with the distance $x_1$ to the leading edge, that is
\begin{equation}
	\Rey_{x_1}=\frac{u_\infty x_1}{\nu}.
\end{equation}
The critical Reynolds number, above which transition to turbulence typically occurs, is approximately $\Rey_{\mathrm{crit}}\approx\num{5E5}$, see \cite{SchlichtingGersten00}. In our described setting, both the free-stream velocity and the domain are fixed and thus only the kinematic viscosity $\nu$ can be used to adjust the Reynolds number. In the following, we set $\nu=\num{2.6E-5}$ as corresponding to $\Rey_{0.5}\approx\num{2E4}$ and thus consider a laminar flow configuration not subject to the onset of turbulence yet still computationally challenging.\\

For the numerical simulations, we always guarantee that SV and TH computations are comparable by using a coarser mesh for SV such that the total number of DOFs is approximately the same. Indeed, Table \ref{tab:BlasiusDOFs} shows that the meshes are chosen such that the SV-FEM is always computationally cheaper. Experience with Taylor--Hood FEM for this test case has revealed that the $\ger{skewTH}_2$ method always yields better results than both EMAC and convective formulation. Thus, we exclusively compare $\ger{skewTH}_2$ and $\ger{SV}_2$. \\

\begin{figure}[h]
\centering
\subfigure{\includegraphics[width=0.975\textwidth]{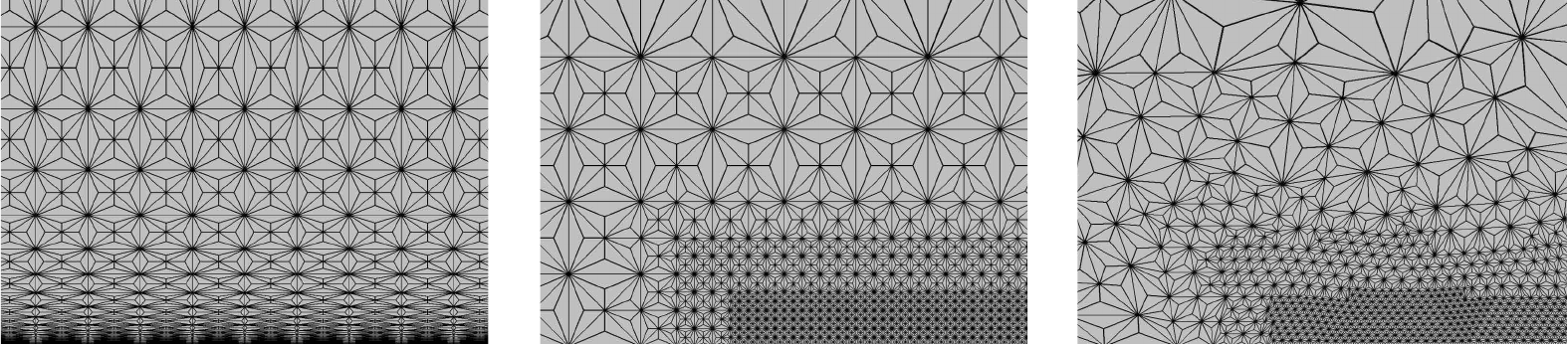}}
\caption{Cutout of barycentre-refined meshes in the range $\rb{-0.28,0.28}\times\rb{0,0.3}$ for the Blasius simulations near the leading edge of the plate. Column-wise: Anisotropic (aspect ratio up to 500), structured isotropic and unstructured isotropic meshes.}
\label{fig:BlasiusMeshes}
\end{figure}

In the following, we compare three different types of meshes: anisotropically refined, structured isotropically refined and unstructured isotropically refined. In Figure \ref{fig:BlasiusMeshes}, cutouts of the meshes used for SV-FEM are shown, where the refinement is always made towards the flat plate, in order to better resolve the processes occurring in the boundary layer. For the TH-FEM, except the barycentre-refinement, the isotropically refined meshes are basically similar. However, the anisotropic mesh for the TH method has a maximum aspect ratio of only 50, whereas the barycentre-refined mesh for the SV method has a maximum aspect ratio of 500. Our experiments with anisotropic meshes revealed that exactly divergence-free methods allow for much more anisotropy (here, ten times higher aspect ratios) without corrupting either convergence or accuracy. \\

\begin{table}[h]
\caption{Number of DOFs arising from different meshes; see Figure \ref{fig:BlasiusMeshes} for barycentre-refined meshes.}
\label{tab:BlasiusDOFs}
\centering 
\begin{tabular}{ccccc} 
\toprule
Method				& Type of mesh 	& Anisotropic	& Structured isotropic 	& Unstructured isotropic\\ 
\otoprule
Taylor--Hood 		& standard		& \num{45503}	& \num{59928} 			& \num{54348}\\ 
Scott--Vogelius 	& barycentre 	& \num{40834}	& \num{53466}			& \num{51954}\\
\bottomrule
\end{tabular}
\end{table}

\begin{figure}[h]
\centering
\subfigure{\includegraphics[width=0.49\textwidth]{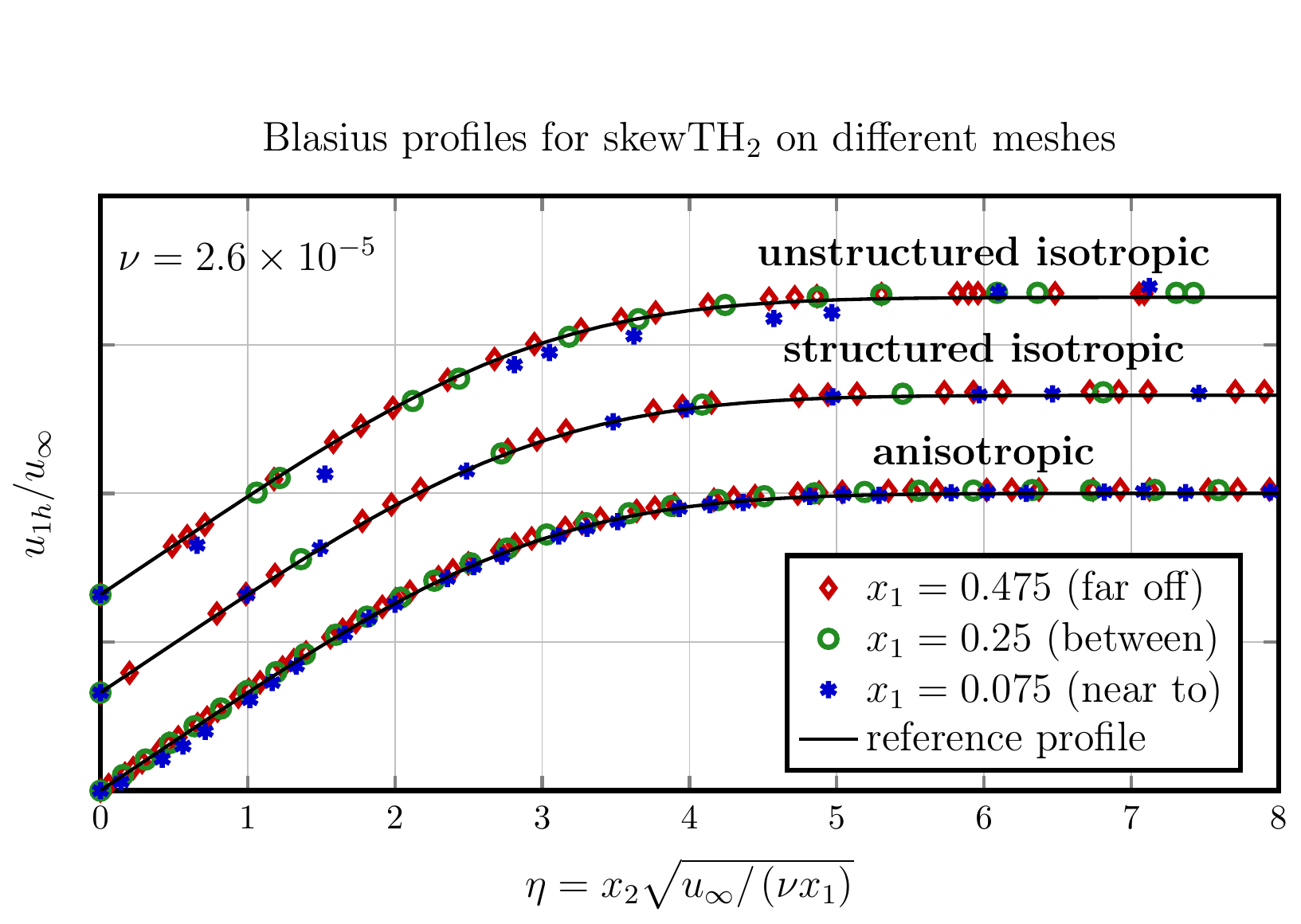}} \goodgap
\subfigure{\includegraphics[width=0.49\textwidth]{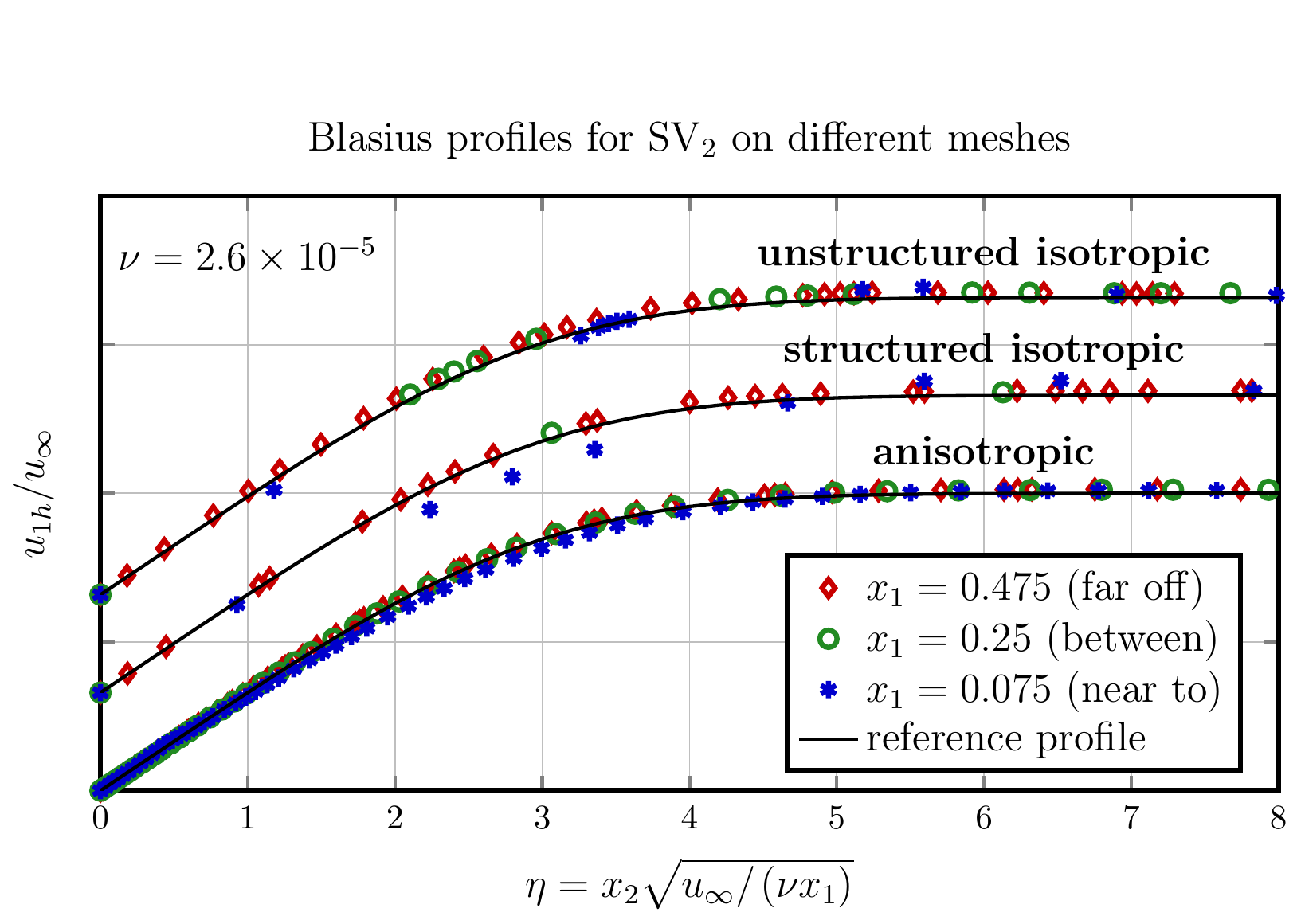}}
\caption{Comparison of computed Blasius profiles ($\ger{skewTH}_2$ and $\ger{SV}_2$) with the self-similar reference solution. Shown are the normalised velocities $u_{1h}/u_\infty$ in $x_1$-direction against the similarity variable $\eta$. For each method, the profiles are evaluated at three different lines $x_1\in\set{0.075,0.25,0.475}$ corresponding to positions far off, in between and near to the leading edge. The results for the different meshes are vertically shifted only for better clarity.}
\label{fig:BlasiusProfiles}
\end{figure}

In Figure \ref{fig:BlasiusProfiles} a comparison of computed Blasius profiles ($\ger{skewTH}_2$ and $\ger{SV}_2$) with the self-similar reference solution on the previously introduced meshes for the normalised  $u_{1h}$-component is shown. For each method, the profiles are evaluated at three different lines $x_1\in\set{0.075,0.25,0.475}$, thereby shedding light on the situation at different distances away from the leading edge of the plate. Let us first focus on apparent similarities between the exactly divergence-free and non divergence-free FEM. First of all, we note that both methods yield satisfactory results for this kind of laminar boundary layer flows. Furthermore, with increasing distance from the leading edge at $\rb{0,0}^\dag$, all numerical solutions become more accurate. Given the nature of the flow which  impinges on the edge of the plate, this is not too surprising. Even though our analysis does not \emph{per se} hold for anisotropic meshes, both $\ger{skewTH}_2$ and $\ger{SV}_2$ yield the most accurate solutions with comparably few DOFs on such meshes. This is surprising because the anisotropic meshes, shown in Figure \ref{fig:BlasiusMeshes}, contain elements with angles very close to $\SI{180}{\degree}$; the barycentre-refined meshes even more so than the standard meshes. \\

The differences between $\ger{skewTH}_2$ and $\ger{SV}_2$ for the Blasius problem are more difficult to recognise. In Figure \ref{fig:BlasiusProfiles} we see that $\ger{skewTH}_2$ yields the worse profiles on the unstructured isotropic mesh, whereas the results from the $\ger{SV}_2$ method on this kind of mesh are satisfying. Vice versa, on the structured isotropic meshes the $\ger{skewTH}_2$ method outperforms $\ger{SV}_2$ slightly. \\

Summarising, we conclude that, despite the rather unfavourable need for barycentre-refined meshes, Scott--Vogelius FEM are at least not inferior to standard Taylor--Hood FEM for this kind of laminar boundary layer flow. For the same amount of computational cost, both methods produce reasonably accurate results on all kinds of considered meshes. However, an advantage of exactly divergence-free FEM over non divergence-free FEM can indeed be observed, since the former allows for stronger anisotropy in the underlying mesh. Thus, we believe that exact mass conservation plays a role in accurately resolving the processes occurring in boundary layers.

\section{Summary and conclusions}	\label{sec:SumAndConcl}

We have considered weakly divergence-free, $\HM{1}$-conforming and inf-sup stable FE approximations of the time-dependent Navier--Stokes problem. Whilst this work has primarily concerned itself with unstabilised Galerkin-FEM, we also pointed to a possible use of edge-stabilisation in the case of high Reynolds number flows. Because we were dealing with pointwise divergence-free approximations, there was no need for any mechanism stabilising the divergence of the velocity approximation (as, for example, in the case of grad-div stabilisation). \\

A stability and convergence analysis was carried out for the continuous-in-time problem, where the focus was on statements concerning the velocity approximation. Similarly to \cite{BurmanFernandez07,ArndtEtAl15}, we have shown that the Gronwall constant is independent of the Reynolds number and therefore our error estimates also hold true verbatim for the incompressible Euler equations. In contrast to these works, our analysis is based on the use of the discrete Stokes projection, which allowed us to fall back on known results for the stationary Stokes problem. Even more to the point, we recognised that Stokes and Ritz projection coincide for discretely divergence-free functions. \\

Assuming that the exact solution fulfils $\uu\in\Lp{\infty}\rb{0,\tend;\WMP{1}{\infty}}\cap\Lp{2}\rb{0,\tend;\HM{k+1}}$, we have derived error estimates of order $\mathcal{O}\rb{h^{k}}$ for the velocity approximation  that are independent of the pressure. Additionally, under the above assumption, we have been able to show that for the discrete velocity, $\uu_h\in\Lp{\infty}\rb{0,\tend;\LP{\infty}}$ holds true. The ability to bound $\uu_h$ in the corresponding norm has been crucial for obtaining $\mathcal{O}\rb{h^{k}}$ error bounds for the discrete pressure, provided $p\in\Lp{2}\rb{0,\tend;\Hm{k}}$. \\ 

We also conducted numerical experiments aimed at comparing the performance of the standard Taylor--Hood method to the Scott--Vogelius method, as it is a representative of the class of exactly divergence-free and conforming FEM. The first two studies were concerned with viscous vortex dynamics, where we took initial velocities that solve the incompressible Euler equations and let them evolve within a viscous incompressible flow. For the third study, the laminar Blasius boundary layer flow has been revisited. \\

The vortex dynamics examples have revealed that, for such flows, divergence-free approximations are superior to standard ones. In particular, we have shown that the conservation of kinetic energy, as well as linear and angular momentum, play a key role for the stability of the particular FEM, and that edge-stabilisation clearly improves the quality of the approximation. Lastly, as a proof of concept, the Blasius problem has revealed that Scott--Vogelius FEM can be used successfully even on highly anisotropic meshes where Taylor--Hood elements fail to yield a converged solution. \\

Finally, let us briefly comment on possible directions for future research. Concerning the error analysis, the question as to whether it is possible to recover half an order in the velocity estimates remains open; however, we found that it is directly connected to the smoothness of $\partial_t\uu$ and the presence of the nonlinear inertia term in the momentum balance. Furthermore, it would be interesting to see whether the reconstruction methods of \cite{LinkeEtAl16,LinkeMerdon16,LedererEtAl16} also allow for semi-robust error estimates; pressure-robustness is guaranteed by construction.

\section*{Acknowledgements}	
We gratefully acknowledge the comments and suggestions about this manuscript from the anonymous reviewers. Furthermore, special thanks go to David Allen Martin II for proofreading and improving the English language used in this work.

\appendix
\renewcommand*{\thesection}{\Alph{section}}
\section{Analysis of the pressure approximation}	 \label{sec:ErrorEstimatesPre}

For the analysis of the pressure approximation, in some places, ideas from \cite{BurmanFernandez07,ArndtEtAl15,DallmannEtAl16,FrutosEtAl16} have been used. Note that, for the sake of brevity, we omit stability estimates for $p_h$. However, the well-posedness of the problem for the pressure can be shown similarly as for the velocity in Section \ref{sec:StabWellPosed}.

\begin{thmLem}[Pressure discretisation error estimate] \label{lem:PresDiscErr}
Under the assumptions of Lemma \ref{lem:uhInfinity}, we obtain the following error estimate for the pressure:
\begin{equation} 
\norm{\pERR}_{\Lp{2}\rb{0,\tend;\Lp{2}}}^2  
	\leqslant \frac{C}{\beta_h^2}\int_0^T \sqb{
	 \norm{\pETA\rb{\tau}}_\Lp{2}^2
	+\norm{\partial_t\uXI\rb{\tau}}_{\VV^*}^2
	+\nu\norm{\uXI\rb{\tau}}_e^2
	+K^2 \norm{\uXI\rb{\tau}}_\LP{2}^2 
	}\dtau
\end{equation}
where $K=\norm{\uu}_{\Lp{\infty}\rb{0,\tend;\LP{\infty}}}+\norm{\uu_h}_{\Lp{\infty}\rb{0,\tend;\LP{\infty}}}$.
\end{thmLem}

\begin{thmProof}
	Contrary to Section \ref{sec:ErrorEstimatesVel}, we now use the $\Lp{2}$-projection as an approximation operator for the pressure:
\begin{align}
\uu-\uu_h&=\rb{\uu-\PIs\uu}+\rb{\PIs\uu-\uu_h}=\uETA + \uERR =\uXI,	\\
p-p_h&=\rb{p-\pi_0 p}+\rb{\pi_0 p -p_h}=\pETA+\pERR=\pXI,
\end{align}	

It is well-known that the discrete inf-sup condition \eqref{eq:InfSup} is equivalent to \cite{GiraultRaviart86}
\begin{equation} \label{eq:VelLift}
	\forall q_h\in\Q_h \quad\exists! \vvl\in\VV_h \quad\text{such that}\quad
	\nabla\cdot\vvl=q_h \quad\text{and}\quad 
	\norm{\nabla\vvl}_\LP{2} \leqslant \beta_h^{-1}\norm{q_h}_\Lp{2}. 	
\end{equation}

Denoting the velocity lifting of the discretisation error $\pERR$ by $\vvl$, we use Cauchy--Schwarz to infer
\begin{subequations}
\begin{align}
	\norm{\pERR}_\Lp{2}^2&=
	\rb{\pERR,\nabla\cdot\vvl}=\rb{\pXI,\nabla\cdot\vvl}-\rb{\pETA,\nabla\cdot\vvl} \leqslant
	\rb{\pXI,\nabla\cdot\vvl} + \norm{\pETA}_\Lp{2}\norm{\nabla\cdot\vvl}_\Lp{2} \\ &\leqslant
	\rb{\pXI,\nabla\cdot\vvl} + \norm{\pETA}_\Lp{2}\norm{\nabla\vvl}_\LP{2}.
\end{align}	
\end{subequations}

Corollary \ref{cor:GalOrtho} with $\rb{\vvl,0}\in\VV_h^\dvg\times\Q_h$ as test functions leads to
\begin{align}
\bra{\partial_t\uXI,\vvl}+a\rb{\uXI,\vvl}+t\rb{\uu;\uu,\vvl}-t\rb{\uu_h;\uu_h,\vvl}+b\rb{\vvl,\pXI} -b\rb{\uXI,0}=0,
\end{align}

from which, using Cauchy--Schwarz for both the duality pairing and $\LTWO$-inner product, we obtain
\begin{align}
	\rb{\pXI,\nabla\cdot\vvl} &\leqslant 
		\norm{\partial_t\uXI}_{\VV^*}\norm{\nabla\vvl}_\LP{2} +
		\nu \norm{\nabla\uXI}_\LP{2}\norm{\nabla\vvl}_\LP{2} +
		t\rb{\uu;\uu,\vvl}-t\rb{\uu_h;\uu_h,\vvl}.
\end{align}	

For the estimation of the difference of convective terms, we add a zero to introduce the error $\uXI$ and then use Lemma \ref{lem:Trilinear} to interchange the second and third argument of the trilinear form: 
\begin{subequations}
\begin{align}
	t\rb{\uu;\uu,\vvl}&-t\rb{\uu_h;\uu_h,\vvl} = 
		\rb{\uu\cdot\nabla\uu,\vvl} - \rb{\uu\cdot\nabla\uu_h,\vvl} + \rb{\uu\cdot\nabla\uu_h,\vvl} - \rb{\uu_h\cdot\nabla\uu_h,\vvl} \\
		&= \rb{\uu\cdot\nabla\sqb{\uu-\uu_h},\vvl} + \rb{\sqb{\uu-\uu_h}\cdot\nabla\uu_h,\vvl} =
		t\rb{\uu;\uXI,\vvl}+t\rb{\uXI;\uu_h,\vvl} \\
		&=-t\rb{\uu;\vvl,\uXI}-t\rb{\uXI;\vvl,\uu_h}
\end{align}	
\end{subequations}

Now, we again invoke Lemma \ref{lem:Trilinear} to estimate this term using
\begin{align} \label{eq:DiffConvShort}
	 \abs{t\rb{\uu;\vvl,\uXI}-t\rb{\uXI;\vvl,\uu_h}}  \leqslant
	\rb{\norm{\uu}_\LP{\infty}+\norm{\uu_h}_\LP{\infty}}\norm{\uXI}_\LP{2}\norm{\nabla\vvl}_\LP{2}.
\end{align}	

Collecting the estimates and using the equivalent formulation of the discrete inf-sup condition yields
\begin{subequations}
\begin{align}
	\beta_h \norm{\pERR}_\Lp{2} &\leqslant
	\frac{\norm{\pERR}_\Lp{2}^2}{\norm{\nabla\vvl}_\LP{2}} = 
	\frac{\rb{\pERR,\nabla\cdot\vvl}}{\norm{\nabla\vvl}_\LP{2}} \\ &\leqslant 
	 \norm{\pETA}_\Lp{2}  + \norm{\partial_t\uXI}_{\VV^*}+
		\nu \norm{\nabla\uXI}_\LP{2} +
		\rb{\norm{\uu}_\LP{\infty}+\norm{\uu_h}_\LP{\infty}}\norm{\uXI}_\LP{2}.
\end{align}	
\end{subequations}

Squaring this inequality and integration over $\rb{0,t}$, for $0\leqslant t\leqslant T$, concludes the proof.

\end{thmProof}

\begin{thmThe}[Pressure discretisation error convergence rate] 
Under the assumptions of the previous lemma, assume a smooth solution according to 
\begin{align}
	\uu&\in\Lp{\infty}\rb{0,\tend;\WMP{1}{\infty}}\cap\Lp{\infty}\rb{0,\tend;\HM{r}}, \quad  
	\partial_t\uu\in\Lp{2}\rb{0,\tend;\HM{r-1}}, \quad
	p\in\Lp{2}\rb{0,\tend;\Hm{s}}.	
\end{align}
Then, with $r_\uu=\min\set{r,k+1}$ and $r_p=\min\set{s,\ell+1}$, we obtain the following convergence rate:
\begin{equation} 
\norm{\pERR}_{\Lp{2}\rb{0,\tend;\Lp{2}}}^2  
	\leqslant \frac{C}{\beta_h^2}\sqb{ h^{2\rb{r_\uu-1}}K^2e^{C_\uu\tend}
				\int_0^\tend \sqb{\abs{\uu\rb{\tau}}_\HM{r_\uu}^2
				+\abs{\partial_t\uu\rb{\tau}}_\HM{r_\uu-1}^2} \dtau
		+ h^{2r_p}\int_0^T \abs{p\rb{\tau}}_\Hm{r_p}^2 \dtau}
\end{equation}
where $K=\norm{\uu}_{\Lp{\infty}\rb{0,\tend;\LP{\infty}}}+\norm{\uu_h}_{\Lp{\infty}\rb{0,\tend;\LP{\infty}}}$. The Gronwall constant $C_\uu$ is given in Theorem \ref{thm:VelDiscErr}.
\end{thmThe}

\begin{thmProof}
The first part of the proof is a consequence of the approximation properties of the $\Ltwo$-projection \eqref{eq:pProj} and Corollary \ref{cor:VelConvRate} (use the triangle inequality to extend the result to the full error $\uu-\uu_h$):
\begin{subequations}
	\begin{align}
		\int_0^T \norm{\pETA\rb{\tau}}_\Lp{2}^2 \dtau
			&\leqslant Ch^{2r_p} \int_0^T \abs{p\rb{\tau}}_\Hm{r_p}^2 \dtau \\
		\nu\int_0^T \norm{\uXI\rb{\tau}}_e^2 \dtau
			&\leqslant C\nu h^{2\rb{r_\uu-1}}e^{C_\uu\tend} 
				\int_0^\tend \sqb{\abs{\uu\rb{\tau}}_\HM{r_\uu}^2
				+\abs{\partial_t\uu\rb{\tau}}_\HM{r_\uu-1}^2} \dtau \\
		K^2\int_0^T \norm{\uXI\rb{\tau}}_\LP{2}^2 \dtau
			&\leqslant CK^2h^{2\rb{r_\uu-1}}e^{C_\uu\tend}
				\int_0^\tend \sqb{\abs{\uu\rb{\tau}}_\HM{r_\uu}^2
				+\abs{\partial_t\uu\rb{\tau}}_\HM{r_\uu-1}^2} \dtau
	\end{align}
\end{subequations}

These are the terms in Lemma \ref{lem:PresDiscErr} that can be estimated directly. The remaining task is then to obtain an estimate of the time derivative of the velocities in the dual norm of $\HM{1}$. First of all, let us concentrate on the dual norm for the discretisation error, that is
\begin{equation}
	\norm{\partial_t\uERR}_{\VV^*}
		= \sup_{\vv\in\VV\backslash\set{\zero}}
			\frac{\abs{\bra{\partial_t\uERR,\vv}}}{\norm{\nabla\vv}_\LP{2}}.
\end{equation}

Introducing the linear operator $A_h\colon\VV_h^\dvg\to\VV_h^\dvg$ defined by
\begin{equation}
	\rb{A_h\vv_h,\ww_h}=\rb{\nabla\vv_h,\nabla\ww_h},\quad\forall\vv_h,\ww_h\in\VV_h^\dvg,
\end{equation}

it is shown in \cite{FrutosEtAl16}, for Oseen's equations, that $\norm{\partial_t\uERR}_{\VV^*}\leqslant C\norm{A_h^{-\half}\partial_t\uERR}_\LP{2}$. For this estimate, an inverse inequality must be used. In the remainder of the proof, we derive a bound for the last norm, thereby extending some ideas from \cite{FrutosEtAl16} to the nonlinear Navier--Stokes case. \\

Based on Corollary \ref{cor:GalOrtho}, we use the same reasoning as in the beginning of the proof of Theorem \ref{thm:VelDiscErr} to infer that for all $\vv_h\in\VV_h^\dvg$,
\begin{align} \label{eq:ErrVhdiv}
	\rb{\partial_t\uERR,\vv_h}+a\rb{\uERR,\vv_h} =
	-\rb{\partial_t\uETA,\vv_h}-\sqb{t\rb{\uu;\uu,\vv_h} - t\rb{\uu_h;\uu_h,\vv_h}}
\end{align}	

holds true. Note the absence of the pressure, thanks to the weakly divergence-free FEM. In the next step, we want to reformulate this equality in terms of an operator equality. To this end, we additionally need the discrete Leray projector $\PI_0^\dvg\colon\LP{2}{\rb{\Omega}}\to \VV_h^\dvg$ ($\LTWO$-projection onto $\VV_h^\dvg$), defined by
\begin{equation}
	\rb{\PI_h^\dvg\vv,\ww_h}=\rb{\vv,\ww_h},\quad\forall\ww_h\in\VV_h^\dvg.
\end{equation}

Then, equation \eqref{eq:ErrVhdiv} can be recast in $\VV_h^\dvg$ as
\begin{equation}
	\partial_t\uERR 
		= 	-\nu A_h\uERR 
			- \PI_0^\dvg\rb{\partial_t\uETA}
			- \PI_0^\dvg\rb{\uu\cdot\nabla\uu-\uu_h\cdot\nabla\uu_h}.
\end{equation}

After applying $A_h^{-\half}$ and taking the $\LTWO$-norm, the triangle inequality and the fact that
\begin{equation}
	\norm{A_h^{-\half}\PI_0^\dvg \gbld}_\LP{2}
		\leqslant \norm{\gbld}_{\VV^*},
		\quad \forall \gbld\in\LP{2}{\rb{\Omega}},
\end{equation}

cf. \cite{FrutosEtAl16}, leads to
\begin{equation}
	\norm{A_h^{-\half}\partial_t\uERR}_\LP{2} 
		\leqslant \nu\norm{A_h^\half\uERR}_\LP{2} 
			+ \norm{\partial_t\uETA}_{\VV^*}
			+ \norm{\rb{\uu\cdot\nabla}\uu-\rb{\uu_h\cdot\nabla}\uu_h}_{\VV^*}.
\end{equation}

Squaring and integrating this inequality in $\sqb{0,t}$, for $0\leqslant t\leqslant \tend$, yields
\begin{subequations}
\begin{align}
	\int_0^t \norm{A_h^\half\partial_t\uERR\rb{\tau}}_\LP{2}^2 \dtau
		\lesssim &\int_0^t \nu^2\norm{A_h^\half\uERR\rb{\tau}}_\LP{2}^2 \dtau
			+ \int_0^t \norm{\partial_t\uETA\rb{\tau}}_{\VV^*}^2 \dtau \\
			+ &\int_0^t \norm{\rb{\uu\cdot\nabla\uu-\uu_h\cdot\nabla\uu_h}\rb{\tau}}_{\VV^*}^2 \dtau.		
\end{align}	
\end{subequations}

Using that $\norm{A_h^\half\vv_h}_\LP{2}=\norm{\nabla\vv_h}_\LP{2}$ for all $\vv_h\in\VV_h^\dvg$, the first two terms on the right-hand side can be handled as in Corollary \ref{cor:VelConvRate}:
\begin{subequations}
\begin{align}
	\int_0^t \nu^2\norm{A_h^\half\uERR\rb{\tau}}_\LP{2}^2 \dtau
		&\leqslant C\nu h^{2\rb{r_\uu-1}}e^{C_\uu\tend} 
				\int_0^\tend \sqb{\abs{\uu\rb{\tau}}_\HM{r_\uu}^2
				+\abs{\partial_t\uu\rb{\tau}}_\HM{r_\uu-1}^2} \dtau \\
	\int_0^t \norm{\partial_t\uETA\rb{\tau}}_{\VV^*}^2 \dtau 
		&\leqslant C \int_0^t \norm{\partial_t\uETA\rb{\tau}}_\LP{2}^2 \dtau
		\leqslant C h^{2\rb{r_\uu-1}}\int_0^t \abs{\partial_t\uu\rb{\tau}}_\HM{r_\uu-1}^2 \dtau
\end{align}	 
\end{subequations}

For the remaining dual norm of the convective terms, we have
\begin{equation}
	\norm{\uu\cdot\nabla\uu-\uu_h\cdot\nabla\uu_h}_{\VV^*}
		= \sup_{\vv\in\VV\backslash\set{\zero}}
			\frac{\abs{t\rb{\uu;\uu,\vv}-t\rb{\uu_h;\uu_h,\vv}}}{\norm{\nabla\vv}_\LP{2}}
		\leqslant \rb{\norm{\uu}_\LP{\infty}+\norm{\uu_h}_\LP{\infty}}\norm{\uXI}_\LP{2},
\end{equation}

which follows from the same computation as \eqref{eq:DiffConvShort}. Finally, again using Corollary \ref{cor:VelConvRate} results in
\begin{subequations}
\begin{align}
	\int_0^t \norm{\rb{\uu\cdot\nabla\uu-\uu_h\cdot\nabla\uu_h}\rb{\tau}}_{\VV^*}^2 \dtau
		&\leqslant K^2\int_0^T \norm{\uXI\rb{\tau}}_\LP{2}^2 \dtau \\
			&\leqslant CK^2h^{2\rb{r_\uu-1}}e^{C_\uu\tend}
				\int_0^\tend \sqb{\abs{\uu\rb{\tau}}_\HM{r_\uu}^2
				+\abs{\partial_t\uu\rb{\tau}}_\HM{r_\uu-1}^2} \dtau.
\end{align}	
\end{subequations}

Combining the above estimates concludes the proof.
\end{thmProof}

\def\bibsection{\section*{\fontsize{10.5}{17}\bfseries References}}

\bibliography{DivFreeFEM4TINS_BibTeX}
\bibliographystyle{alphaabbr}

\end{document}